\documentclass[a4paper]{article}%
\usepackage{graphicx}
\usepackage{geometry}
\usepackage{amsmath}
\usepackage{amssymb}
\usepackage{amsfonts}
\usepackage{amsthm,amscd}%
\def\figurename{Figure}
\makeatletter
\renewcommand{\fnum@figure}[1]{\figurename~\thefigure.}
\makeatother
\def\tablename{Table}
\makeatletter
\renewcommand{\fnum@table}[1]{\tablename~\thetable.}
\makeatother
\def \bop {\noindent\textbf{Proof. }}
\def \eop {\hbox{}\nobreak\hfill
\vrule width 2mm height 2mm depth 0mm
\par \goodbreak \smallskip}
\newtheorem{theorem}{Theorem}[section]
\newtheorem{lemma}[theorem]{Lemma}
\newtheorem{corollary}[theorem]{Corollary}

\theoremstyle{definition}
\newtheorem{definition}[theorem]{Definition}

\theoremstyle{remark}
\newtheorem{remark}[theorem]{Remark}
\numberwithin{equation}{section}

\begin{document}

\title{The Maximum Principle in Time-Inconsistent LQ Optimal Control Problem for Jump Diffusions}
\author{\textbf{Ishak Alia}{ }\thanks{ Laboratory of Applied Mathematics, University
Mohamed Khider, Po. Box 145 Biskra (07000), Algeria. E-mail:
ishak.alia@hotmail.com}
\and \textbf{Farid Chighoub}{ }\thanks{ Laboratory of Applied Mathematics,
University Mohamed Khider, Po. Box 145 Biskra (07000), Algeria. E-mail:
f.chighoub@univ-biskra.dz, chighoub\_farid@yahoo.fr}
\and \textbf{Ayesha Sohail }\thanks{Department of Mathematics, Comsats Institute of
Information Technology, Lahore, Pakistan. E-mail address:
asohail@ciitlahore.edu.pk}}
\maketitle

\begin{abstract}
In this paper, we consider a general time-inconsistent optimal control problem
for a non homogeneous linear system, in which its state evolves according to a
stochastic differential equation with deterministic coefficients, when the
noise is driven by a Brownian motion and an independent Poisson point process.
The running and the terminal costs in the objective functional, are explicitly
dependent on some general discounting coefficients which cover the
non-exponential and the hyperbolic discounting situations. Furthermore, the
presence of some quadratic terms of the conditional expectation of the state
process as well as a state-dependent term in the objective functional makes
the problem time-inconsistent. Open-loop Nash equilibrium controls are
constructed instead of optimal controls, by using a version of the stochastic
maximum principle approach. This approach involves a stochastic system that
consists of a flow of forward-backward stochastic differential equations and
an equilibrium condition. As an application, we study some concrete examples.

\end{abstract}

\textbf{Keys words}: Stochastic maximum principle, Time inconsistency, Linear
quadratic control problem, Equilibrium control, Variational inequality.


\section{Introduction}

Time-inconsistent stochastic control problems have received remarkable
attention in the recent years. The risk aversion attitude of a mean-variance
investor \cite{8}, \cite{13} and \cite{18}, such as the portfolio optimization
with non-exponential discount function \cite{1} and \cite{2}, provide two
well-known examples of time-inconsistency in mathematical finance. Motivated
by these practical examples, this paper studies optimality conditions for
time-inconsistent linear quadratic stochastic control problem, where the state
is described by a n-dimensional non homogeneous controlled SDE with jump
processes, defined on a complete filtered probability space $(\Omega
,\mathcal{F},\left(  \mathcal{F}_{t}\right)  _{t\in\left[  0,T\right]
},\mathbb{P})$%
\begin{equation}
\left\{
\begin{array}
[c]{l}%
dX\left(  s\right)  =\left\{  A\left(  s\right)  X\left(  s\right)  +B\left(
s\right)  u\left(  s\right)  +b\left(  s\right)  \right\}  ds+%
{\displaystyle\sum\limits_{j=1}^{d}}
\left\{  C_{j}\left(  s\right)  X\left(  s\right)  +D_{j}\left(  s\right)
u\left(  s\right)  +\sigma_{j}\left(  s\right)  \right\}  dW^{j}\left(
s\right) \\
\text{ \ \ \ \ \ \ \ \ \ \ }+%
{\displaystyle\int\limits_{Z}}
\left\{  E\left(  s,z\right)  X\left(  s-\right)  +F\left(  s,z\right)
u\left(  s\right)  +c\left(  s,z\right)  \right\}  \tilde{N}\left(
ds,dz\right)  ,\text{\ }s\in\left[  0,T\right]  ,\\
X\left(  0\right)  =x_{0}\left(  \in%
\mathbb{R}
^{n}\right)  \text{.}%
\end{array}
\right.  \tag{1.1}%
\end{equation}

The coefficients $A\left(  .\right)  ,B\left(  .\right)  ,b\left(  .\right)
,C_{j}\left(  .\right)  ,D_{j}\left(  .\right)  ,\sigma_{j}\left(  .\right)
,E\left(  .,.\right)  ,F\left(  .,.\right)  $ and $c\left(  .,.\right)  $ are
deterministic matrix-valued functions of suitable sizes. As time evolves, it
is natural to consider the linear controlled stochastic differential equation
starting from the situation $\left(  t,\xi\right)  \in\left[  0,T\right]
\times\mathbb{L}^{2}\left(  \Omega,\mathcal{F}_{t},\mathbb{P};%
\mathbb{R}
^{n}\right)  $%
\begin{equation}
\left\{
\begin{array}
[c]{l}%
dX\left(  s\right)  =\left\{  A\left(  s\right)  X\left(  s\right)  +B\left(
s\right)  u\left(  s\right)  +b\left(  s\right)  \right\}  ds+%
{\displaystyle\sum\limits_{j=1}^{d}}
\left\{  C_{j}\left(  s\right)  X\left(  s\right)  +D_{j}\left(  t\right)
u\left(  s\right)  +\sigma_{j}\left(  s\right)  \right\}  dW^{j}\left(
s\right) \\
\text{ \ \ \ \ \ \ \ \ \ \ }+%
{\displaystyle\int\limits_{Z}}
\left\{  E\left(  s,z\right)  X\left(  s-\right)  +F\left(  s,z\right)
u\left(  s\right)  +c\left(  s,z\right)  \right\}  \tilde{N}\left(
ds,dz\right)  ,\text{\ }s\in\left[  t,T\right]  ,\\
X\left(  t\right)  =\xi\text{.}%
\end{array}
\right.  \tag{1.2}%
\end{equation}

Under some conditions, for any initial situation $\left(  t,\xi\right)  $ and
any admissible control $u\left(  .\right)  \ $the state equation is uniquely
solvable, we denote by $X\left(  .\right)  =X^{t,\xi}\left(  .;u\left(
.\right)  \right)  $ its solution$,$ for $s\in\left[  t,T\right]  .$ Different
controls $u\left(  .\right)  $ will lead to different solutions $X\left(
.\right)  .$ To measure the performance of $u\left(  .\right)  ,\ $we
introduce the following cost functional%
\begin{equation}%
\begin{array}
[c]{l}%
J\left(  t,\xi,u\left(  .\right)  \right) \\
=\mathbb{E}^{t}\left[
{\displaystyle\int_{t}^{T}}
\dfrac{1}{2}\left(  \left\langle Q\left(  t,s\right)  X\left(  s\right)
,X\left(  s\right)  \right\rangle +\left\langle \bar{Q}\left(  t,s\right)
\mathbb{E}^{t}\left[  X\left(  s\right)  \right]  ,\mathbb{E}^{t}\left[
X\left(  s\right)  \right]  \right\rangle +\left\langle R\left(  t,s\right)
u\left(  s\right)  ,u\left(  s\right)  \right\rangle \right)  ds\right. \\
\text{ \ \ \ \ \ \ \ \ \ \ \ }+\left\langle \mu_{1}\left(  t\right)  \xi
+\mu_{2}\left(  t\right)  ,X\left(  T\right)  \right\rangle \\
\text{ \ \ \ \ \ \ \ \ \ \ \ }+\left.  \dfrac{1}{2}\left(  \left\langle
G\left(  t\right)  X\left(  T\right)  ,X\left(  T\right)  \right\rangle
+\left\langle \bar{G}\left(  t\right)  \mathbb{E}^{t}\left[  X\left(
T\right)  \right]  ,\mathbb{E}^{t}\left[  X\left(  T\right)  \right]
\right\rangle \right)  \right]  .
\end{array}
\tag{1.3}%
\end{equation}

The coefficients $Q\left(  .,.\right)  ,\bar{Q}\left(  .,.\right)  ,R\left(
.,.\right)  ,G\left(  .\right)  ,\bar{G}\left(  .\right)  ,\mu_{1}\left(
.\right)  $ and $\mu_{2}\left(  .\right)  $ are deterministic matrix-valued
functions of suitable sizes, which explicitly depend on the initial time $t$
in some general way. Our objective in this paper, is to investigate a general
discounting linear quadratic optimal control problem for jump diffusions,
which is time-inconsistent in the sense that, it does not satisfy the Bellman
optimality principle, since a restriction of an optimal control for a specific
initial pair on a later time interval might not be optimal for that
corresponding initial pair. The novelty of this work lies in the fact that,
our calculations are not limited to the exponential discounting framework, the
time-inconsistency of the LQ optimal controls that we are going to consider,
is due to the presence of some general discounting coefficients, involving the
so-called hyperbolic discounting situations. In addition, the presence of some
quadratic terms of the expected controlled state process, in either the
running cost or the terminal cost, make the problem time-inconsistent, this
can be motivated by the reward term in the mean-variance portfolio choice
model. The term $\mu_{1}\left(  t\right)  \xi+\mu_{2}\left(  t\right)  $ stems
from a state-dependent utility function in economics \cite{18}$.$ Each of
these terms introduces time-inconsistency of the underlying model, in somewhat
different ways.

The main difficulty when facing a time-inconsistent optimal control problem is
that, we cannot use the dynamic programming and the standard HJB techniques,
in general. However, the main approach to handle the time-inconsistent optimal
control problems, is by viewing them within a game theoretic framework. Nash
equilibriums are therefore considered instead of optimal solutions, see e.g.
\cite{8}, \cite{19}, \cite{14}, \cite{1},$\ \cite{2},$ \cite{18}, \cite{5},
\cite{10}, \cite{7}, \cite{6}, $\cite{3},$ \cite{20} and \cite{4}. The
fundamental idea is that the control action that the controller makes at every
instant of time, is considered as a game against all the control actions that
the future incarnations of the controller are going to make. Strotz
\cite{3}$,$ was the first who used this game perspective to handle the dynamic
time-inconsistent decision problem on the deterministic Ramsay problem
\cite{6}$.$ Then by capturing the idea of non-commitment, by letting the
commitment period being infinitesimally small, he characterized a Nash
equilibrium strategy. Further work which extend \cite{6}$,$ are \cite{5},
\cite{6}, \cite{7} and \cite{11}$.$ Ekland and Lazrak \cite{1} and Ekland and
Pirvu \cite{2} apply this game perspective to investigate the optimal
investment-consumption problem under general discount functions, in both,
deterministic and stochastic framework. Then, by means of the so-called "local
spike variation" they provide a formal definition of feedback Nash equilibrium
controls in continuous time. The work \cite{8} extends the idea to the
stochastic framework where the controlled process is Markovian. In addition,
an extended HJB equation is derived, along with a verification argument that
characterizes a Markov subgame perfect Nash equilibium. In $\cite{3},$ Yong
studied\ a time-inconsistent deterministic linear quadratic model, and he
derive a closed-loop equilibrium strategie, via a forward ordinary
differential equation coupled with a backward Riccati--Volterra integral
equation. Hu et al \cite{18} studied a time-inconsistent stochastic
linear--quadratic control model, which is originated from the mean-variance
portfolio selection problem with state-dependent risk aversion, and by means
of variational method they derive a general sufficient condition for
equilibria, through a new class of forward-backward stochastic differential
equation (FBSDE in short) along with some equilibrium conditions. In \cite{4}
Yong investigate a time-inconsistent stochastic problem for stochastic
differential equation. By introducing a family of N-person non-cooperative
differential games he characterize a closed-loop equilibrium strategie.

The purpose of this paper is to characterize Nash equilibrium controls for a
general time-inconsistent stochastic linear quadratic optimal control problem.
The objective functional includes the cases of hyperbolic discounting, as well
as, the continuous-time Markowitz's mean-variance portfolio selection problem,
with state-dependent risk aversion. We accentuate that, our model\textbf{
}covers some class of time-inconsistent stochastic LQ optimal control problem
studied by \cite{18}, and some relevant cases appeared in \cite{20}$.$ Note
that, in \cite{18} the weighting matrices do not depend on current time $t$
and in \cite{20} the terminal cost do not depend on current state $\xi.$
Moreover, we have defined the equilibrium controls in open-loop sense (in a
manner similar to \cite{18}), which is different from the feedback form (see
\cite{8}, \cite{13}, \cite{19}, \cite{1}, \cite{2}, \cite{3}, \cite{4} and
\cite{22})$.$

The rest of the paper is organized as follows. In Section 2,\ we describe the
model and formulate the objective. In\ Section 3 we present the first main
result of this work (Theorem 3.2), which characterizes the equilibrium control
via a stochastic system, which involves a flow of forward-backward stochastic
differential equation with jumps (FBSDEJ in short), along with some
equilibrium conditions. In Section 4, by decoupling the flow of the FBSDEJ, we
investigate a feedback representation of the equilibrium control, via some
class of ordinary differential equations, which do not have a symmetry
structure.\ Section 5 is devoted to some applications, we solve a continuous
time mean--variance portfolio selection model and some one-dimensional general
discounting LQ problems. The paper ends with Appendix containing some proofs.

\section{Problem setting}

Let $(\Omega,\mathcal{F},\left(  \mathcal{F}_{t}\right)  _{t\in\left[
0,T\right]  },\mathbb{P})$ be a filtered probability space such that
$\mathcal{F}_{0}$ contains all $\mathbb{P}$-null sets, $\mathcal{F}%
_{T}=\mathcal{F}$ for an arbitrarily fixed finite time horizon $T>0,$ and
$\left(  \mathcal{F}_{t}\right)  _{t\in\left[  0,T\right]  }$ satisfies the
usual conditions. We assume that $\left(  \mathcal{F}_{t}\right)
_{t\in\left[  0,T\right]  }$ is generated by a d-dimensional standard Browian
motion $\left(  W\left(  t\right)  \right)  _{t\in\left[  0,T\right]  }$ and
an independent Poisson measure $N$ on $\left[  0,T\right]  \times Z$ where $Z$
$\subset%
\mathbb{R}
-\left\{  0\right\}  $. We assume that the compensator of $N$ has the form
$\mu\left(  dt,dz\right)  =\theta\left(  dz\right)  dt$ for some positive and
$\sigma-$finite Levy measure on $Z$, endowed with it's Borel $\sigma-$field
$\mathcal{B}\left(  Z\right)  $. We suppose that $%
{\textstyle\int\limits_{Z}}
1\wedge\left\vert z\right\vert ^{2}\theta\left(  dz\right)  <\infty$ and write
$\tilde{N}\left(  dt,dz\right)  =N\left(  dt,dz\right)  -\theta\left(
dz\right)  dt$ for the compensated jump martingal rondom measure of $N.$
Obviously, we have%
\[
\mathcal{F}_{t}=\sigma\left[
{\textstyle\int}
{\textstyle\int_{A\times\left(  0,s\right]  }}
N\left(  dr,de\right)  ;s\leq t,A\in\mathcal{B}\left(  Z\right)  \right]
\vee\sigma\left[  B_{s};s\leq t\right]  \vee\mathcal{N},
\]
where $\mathcal{N}$ denotes the totality of $\nu-$null sets, and $\sigma
_{1}\vee\sigma_{2}$\ denotes the $\sigma-$field generated by $\sigma_{1}%
\cup\sigma_{2}.$

\subsection{Notations}

Throughout this paper, we use the following notations:

\begin{enumerate}
\item[$\bullet$] $S^{n}:$ the set of $n\times n$ symmetric real matrices$.$

\item[$\bullet$] $C^{\top}:$ the transpose of the vector (or matrix) $C.$

\item[$\bullet$] $\left\langle .,.\right\rangle :$ the inner product in some
Euclidean space.
\end{enumerate}

For any Euclidean space $H=%
\mathbb{R}
^{n}$, $%
\mathbb{R}
^{n\times m}$ or $S^{n}$ with Frobenius norm $\left\vert .\right\vert $ we let
for any $t\in\left[  0,T\right]  $

\begin{enumerate}
\item[$\bullet$] $\mathbb{L}^{p}\left(  \Omega,\mathcal{F}_{t},\mathbb{P}%
;H\right)  :=\left\{  \xi:\Omega\rightarrow H\text{ }|\text{ }\xi\text{ is
}\mathcal{F}_{t}-\text{measurable},\text{ with }\mathbb{E}\left[  \left\vert
\xi\right\vert ^{p}\right]  <\infty\right\}  $, for any $p\geq1.$

\item[$\bullet$] $\mathbb{L}^{2}\left(  Z,\mathcal{B}\left(  Z\right)
,\theta;H\right)  :=\left\{  r\left(  .\right)  :Z\rightarrow H\text{ }|\text{
}r\left(  .\right)  \text{ is }\mathcal{B}\left(  Z\right)  -\text{measurable}%
,\text{ with }%
{\textstyle\int_{Z}}
\left\vert r\left(  z\right)  \right\vert ^{2}\theta\left(  dz\right)
<\infty\right\}  .$

\item[$\bullet$] $\mathcal{S}_{%
\mathcal{F}%
}^{2}\left(  t,T;H\right)  :=\left\{  X\left(  .\right)  :\left[  t,T\right]
\times\Omega\rightarrow H\text{ }|\text{ }X\left(  .\right)  \text{ is
}\left(  \mathcal{F}_{s}\right)  _{s\in\left[  t,T\right]  }-\text{adapted}%
,\right.  $
\end{enumerate}

$\ \ \ \ \ \ \ \ \ \ \ \ \ \ \ \ \ \ \ \ \ \ \ \ \ \ \ \ \ \ \ \left.
s\mapsto X(s)\text{ is c\`{a}dl\`{a}g},\text{ with }\mathbb{E}\sup
\limits_{s\in\left[  t,T\right]  }\left\vert X\left(  s\right)  \right\vert
^{2}ds<\infty\right\}  .$

\begin{enumerate}
\item[$\bullet$] $\mathcal{L}_{%
\mathcal{F}%
}^{2}\left(  t,T;H\right)  :=\left\{  X\left(  .\right)  :\left[  t,T\right]
\times\Omega\rightarrow H|X\left(  .\right)  \text{ is }\left(  \mathcal{F}%
_{s}\right)  _{s\in\left[  t,T\right]  }-\text{adapted},\text{ with
}\mathbb{E}\left[
{\textstyle\int_{t}^{T}}
\left\vert X\left(  s\right)  \right\vert ^{2}ds\right]  <\infty\right\}  .$

\item[$\bullet$] $\mathcal{L}_{%
\mathcal{F}%
}^{\theta,2}\left(  \left[  t,T\right]  \times Z;H\right)  :=\left\{  R\left(
.,.\right)  :\left[  t,T\right]  \times\Omega\times Z\rightarrow H\text{
}|\text{ }R\left(  .\right)  \text{ is }\left(  \mathcal{F}_{s}\right)
_{s\in\left[  t,T\right]  }-\text{adapted measurable}\right.  $

$\ \ \ \ \ \ \ \ \ \ \ \ \ \ \ \ \ \ \ \ \ \ \ \ \ \ \ \ \ \ \ \ \ \ \ \ \ \ $%
process on $\left[  t,T\right]  \times\Omega\times Z,\ $with $\left.
\mathbb{E}\left[
{\textstyle\int_{t}^{T}}
{\textstyle\int_{Z}}
\left\vert R\left(  s,z\right)  \right\vert ^{2}\theta\left(  dz\right)
ds\right]  <\infty\right\}  .$

\item[$\bullet$] $C\left(  \left[  0,T\right]  ;H\right)  :=\left\{  f:\left[
0,T\right]  \rightarrow H|\text{ }f\left(  .\right)  \text{ }is\text{
continuous}\right\}  .$

\item[$\bullet$] $D\left[  0,T\right]  :=\left\{  \left(  t,s\right)
\in\left[  0,T\right]  \times\left[  0,T\right]  ,\text{ such that }s\geq
t\right\}  .$

\item[$\bullet$] $C\left(  D\left[  0,T\right]  ;H\right)  :=\left\{  f\left(
.,.\right)  :D\left[  0,T\right]  \rightarrow H\left(  t,s\right)  |\text{
}f\left(  .,.\right)  \text{ }is\text{ continuous}\right\}  .$

\item[$\bullet$] $C^{0,1}\left(  D\left[  0,T\right]  ;H\right)  :=\left\{
f\left(  .,.\right)  :D\left[  0,T\right]  \rightarrow H|\text{ }f\left(
.,.\right)  \text{ and }\dfrac{\partial f}{\partial s}\left(  .,.\right)
\text{\ are continuous}\right\}  .$
\end{enumerate}

\subsection{Problem statement}

We consider a n-dimensional non homogeneous linear controlled jump diffusion
system%
\begin{equation}
\left\{
\begin{array}
[c]{l}%
dX\left(  s\right)  =\left\{  A\left(  s\right)  X\left(  s\right)  +B\left(
s\right)  u\left(  s\right)  +b\left(  s\right)  \right\}  ds+%
{\displaystyle\sum\limits_{j=1}^{d}}
\left\{  C_{j}\left(  s\right)  X\left(  s\right)  +D_{j}\left(  t\right)
u\left(  s\right)  +\sigma_{j}\left(  s\right)  \right\}  dW^{j}\left(
s\right) \\
\text{ \ \ \ \ \ \ \ \ \ \ }+%
{\displaystyle\int\limits_{Z}}
\left\{  E\left(  s,z\right)  X\left(  s-\right)  +F\left(  s,z\right)
u\left(  s\right)  +c\left(  s,z\right)  \right\}  \tilde{N}\left(
ds,dz\right)  ,\text{\ }s\in\left[  t,T\right]  ,\\
X\left(  t\right)  =\xi\text{.}%
\end{array}
\right.  \tag{2.1}%
\end{equation}
where $\left(  t,\xi,u\left(  .\right)  \right)  \in\left[  0,T\right]
\times\mathbb{L}^{2}\left(  \Omega,\mathcal{F}_{t},\mathbb{P};%
\mathbb{R}
^{n}\right)  \times\mathcal{L}_{%
\mathcal{F}%
}^{2}\left(  t,T;%
\mathbb{R}
^{m}\right)  .$ Note that $\mathcal{L}_{%
\mathcal{F}%
}^{2}\left(  t,T;%
\mathbb{R}
^{m}\right)  $ is the space of all admissible strategies. Our aim is to
minimize the following expected discounted cost functional%
\begin{equation}%
\begin{array}
[c]{l}%
J\left(  t,\xi,u\left(  .\right)  \right) \\
=\mathbb{E}^{t}\left[
{\displaystyle\int_{t}^{T}}
\dfrac{1}{2}\left(  \left\langle Q\left(  t,s\right)  X\left(  s\right)
,X\left(  s\right)  \right\rangle +\left\langle \bar{Q}\left(  t,s\right)
\mathbb{E}^{t}\left[  X\left(  s\right)  \right]  ,\mathbb{E}^{t}\left[
X\left(  s\right)  \right]  \right\rangle +\left\langle R\left(  t,s\right)
u\left(  s\right)  ,u\left(  s\right)  \right\rangle \right)  ds\right. \\
\text{ \ \ \ \ \ \ \ \ \ \ \ }+\left\langle \mu_{1}\left(  t\right)  \xi
+\mu_{2}\left(  t\right)  ,X\left(  T\right)  \right\rangle \\
\text{ \ \ \ \ \ \ \ \ \ \ \ }+\left.  \dfrac{1}{2}\left(  \left\langle
G\left(  t\right)  X\left(  T\right)  ,X\left(  T\right)  \right\rangle
+\left\langle \bar{G}\left(  t\right)  \mathbb{E}^{t}\left[  X\left(
T\right)  \right]  ,\mathbb{E}^{t}\left[  X\left(  T\right)  \right]
\right\rangle \right)  \right]  ,
\end{array}
\tag{2.2}%
\end{equation}
over $u\left(  .\right)  \in\mathcal{L}_{%
\mathcal{F}%
}^{2}\left(  t,T;%
\mathbb{R}
^{m}\right)  ,$where $X\left(  .\right)  =X^{t,\xi}\left(  .;u\left(
.\right)  \right)  $ and $\mathbb{E}^{t}\left[  .\right]  =\mathbb{E}\left[
.\left\vert \mathcal{F}_{t}\right.  \right]  .$

We need to impose the following assumptions about the coefficients

\begin{enumerate}
\item[\textbf{(H1)}] \textbf{ }The functions\textbf{\ }$A\left(  .\right)
,C_{j}\left(  .\right)  :\left[  0,T\right]  \rightarrow%
\mathbb{R}
^{n\times n},$ $B\left(  .\right)  ,D_{j}\left(  .\right)  :\left[
0,T\right]  \rightarrow%
\mathbb{R}
^{n\times m},$ $b\left(  .\right)  ,\sigma_{j}\left(  .\right)  :\left[
0,T\right]  \rightarrow%
\mathbb{R}
^{n},$ $E\left(  .,.\right)  :\left[  0,T\right]  \times Z\rightarrow%
\mathbb{R}
^{n\times n},$ $F\left(  .,.\right)  :\left[  0,T\right]  \times Z\rightarrow%
\mathbb{R}
^{n\times m}$, and $c\left(  .,.\right)  :\left[  0,T\right]  \times
Z\rightarrow%
\mathbb{R}
^{n}$ are continuous, and the coefficients on the cost functional\textbf{\ }%
satisfy%
\[
\left\{
\begin{array}
[c]{l}%
Q\left(  .,.\right)  ,\bar{Q}\left(  .,.\right)  \in C\left(  D\left[
0,T\right]  ;S^{n}\right)  ,\medskip\\
R\left(  .,.\right)  \in C\left(  D\left[  0,T\right]  ;S^{m}\right)
,\medskip\\
G\left(  .\right)  ,\bar{G}\left(  .\right)  \in C\left(  \left[  0,T\right]
;S^{n}\right)  ,\medskip\\
\mu_{1}\left(  .\right)  \in C\left(  \left[  0,T\right]  ;%
\mathbb{R}
^{n\times n}\right)  ,\medskip\\
\mu_{2}\left(  .\right)  \in C\left(  \left[  0,T\right]  ;%
\mathbb{R}
^{n}\right)  .
\end{array}
\right.
\]

\item[\textbf{(H2)}] \textbf{ }The functions $R\left(  .,.\right)  ,Q\left(
.,.\right)  $ and $G\left(  .\right)  $\textbf{\ }satisfy\textbf{ }%
\[
R\left(  t,t\right)  \geq0,\text{ }G\left(  t\right)  \geq0,\text{ }\forall
t\in\left[  0,T\right]  ,\text{and }Q\left(  t,s\right)  \geq0,\text{ }%
\forall\left(  t,s\right)  \in D\left[  0,T\right]  .
\]

\end{enumerate}

Under \textbf{(H1) }for any $\left(  t,\xi,u\left(  .\right)  \right)
\in\left[  0,T\right]  \times\mathbb{L}^{2}\left(  \Omega,\mathcal{F}%
_{t},\mathbb{P};%
\mathbb{R}
^{n}\right)  \times\mathcal{L}_{%
\mathcal{F}%
}^{2}\left(  t,T;%
\mathbb{R}
^{m}\right)  ,$ the state equation $\left(  2.1\right)  $ has a unique
solution $X\left(  .\right)  \in\mathcal{S}_{%
\mathcal{F}%
}^{2}\left(  t,T;%
\mathbb{R}
^{n}\right)  ,$ see for example \cite{23}. Moreover, we have the following
estimate%
\[
\mathbb{E}\left[  \sup_{t\leq s\leq T}\left\vert X\left(  s\right)
\right\vert ^{2}\right]  \leq K\left(  1+\mathbb{E}\left[  \left\vert
\xi\right\vert ^{2}\right]  \right)  ,
\]
for some positif constant $K$. The optimal control problem can be formulated
as follows.

\textbf{Problem (LQJ).} \textit{For any given initial pair }$\left(
t,\xi\right)  \in\left[  0,T\right]  \times\mathbb{L}^{2}\left(
\Omega,\mathcal{F}_{t},\mathbb{P};%
\mathbb{R}
^{n}\right)  $\textit{, find a control }$\hat{u}\left(  .\right)
\in\mathcal{L}_{%
\mathcal{F}%
}^{2}\left(  t,T;%
\mathbb{R}
^{m}\right)  $\textit{ such that}%
\[
J\left(  t,\xi,\hat{u}\left(  .\right)  \right)  =\min_{u\left(  .\right)
\in\mathcal{L}_{%
\mathcal{F}%
}^{2}\left(  t,T;%
\mathbb{R}
^{m}\right)  }J\left(  t,\xi,u\left(  .\right)  \right)
\]

\begin{remark}
1) The dependence of the weighting matrices of the current time $t,$ the
dependence of the terminal cost on the current state $\xi$ and the presence of
quadratic terms of the expected controlled state process in the cost
functional make the Problem (LQJ) time-inconsistent.

2) One way to get around the time-inconsistency issue is to consider only
precommitted controls (i.e., the controls are optimal only when viewed at the
initial time).
\end{remark}

\subsection{An example of time-inconsistent optimal control problem}

We present a simple illustration of stochastic optimal control problem which
is time-inconsistent. Our aim is to show that the classical SMP approach is
not efficient in the study of this problem if it's viewed as time-consistent.
For $n=d=1,$ consider the following controlled SDE starting from $\left(
t,x\right)  \in\left[  0,T\right]  \times%
\mathbb{R}
$%
\begin{equation}
\left\{
\begin{array}
[c]{l}%
dX^{t,x}\left(  s\right)  =bu\left(  s\right)  ds+\sigma dW\left(  s\right)
,\text{ }s\in\left[  t,T\right]  ,\\
X^{t,x}\left(  t\right)  =x,
\end{array}
\right.  \tag{2.3}%
\end{equation}
where $b$ and $\sigma$ are real constants. The cost functional given by%
\begin{equation}
J\left(  t,x,u\left(  .\right)  \right)  =\frac{1}{2}\mathbb{E}\left[
\int_{t}^{T}\left\vert u\left(  s\right)  \right\vert ^{2}ds+h\left(
t\right)  \left(  X^{t,x}\left(  T\right)  -x\right)  ^{2}\right]  , \tag{2.4}%
\end{equation}
where $h\left(  .\right)  :\left[  0,T\right]  \rightarrow\left(
0,\infty\right)  ,$ is a general deterministic non-exponential discount
function satisfying $h\left(  0\right)  =1,$ $h\left(  s\right)  \geq0$ and
$\int_{0}^{T}h\left(  t\right)  dt<\infty$. We want to address the following
stochastic control problem.

\textbf{Problem (E).} \textit{For any given initial pair }$\left(  t,x\right)
\in\left[  0,T\right]  \times%
\mathbb{R}
$\textit{, find a control }$\bar{u}\left(  .\right)  \in\mathcal{L}_{%
\mathcal{F}%
}^{2}\left(  t,T;%
\mathbb{R}
\right)  $\textit{ such that }%
\[
J\left(  t,x,\bar{u}\left(  .\right)  \right)  =\min_{u\left(  .\right)
\in\mathcal{L}_{%
\mathcal{F}%
}^{2}\left(  t,T;%
\mathbb{R}
\right)  }J\left(  t,x,u\left(  .\right)  \right)  ,
\]

At a first stage, we consider Problem (E) as a standard time consistent
stochastic linear quadratic problem. Since $J\left(  t,x,.\right)  $ is convex
and coercive, there exists then a unique optimal control for this problem for
each fixed initial pair $\left(  t,x\right)  \in\left[  0,T\right]  \times%
\mathbb{R}
.$ Notice that the usual Hamiltonian associated to this problem is
$\mathbb{H}:\left[  0,T\right]  \times%
\mathbb{R}
^{4}\rightarrow%
\mathbb{R}
$ such that for every $\left(  s,y,v,p,q\right)  \in\left[  0,T\right]  \times%
\mathbb{R}
^{4}$ we have%
\[
\mathbb{H}\left(  s,y,v,p,q\right)  =pbv+\sigma q-\frac{1}{2}v^{2},
\]

Let $u^{t,x}\left(  .\right)  $ be an admissible control for $\left(
t,x\right)  \in\left[  0,T\right]  \times%
\mathbb{R}
.$ Then the corresponding first order and second order adjoint equations are
given respectively by%
\[
\left\{
\begin{array}
[c]{l}%
dp^{t,x}\left(  s\right)  =q^{t,x}\left(  s\right)  dW\left(  s\right)
,\text{ }s\in\left[  t,T\right]  ,\\
p^{t,x}\left(  T\right)  =-h\left(  t\right)  \left(  X^{t,x}\left(  T\right)
-x\right)  ,
\end{array}
\right.
\]
and
\[
\left\{
\begin{array}
[c]{l}%
dP^{t,x}\left(  s\right)  =Q^{t,x}\left(  s\right)  dW\left(  s\right)
,\text{ }s\in\left[  t,T\right]  ,\\
P^{t,x}\left(  T\right)  =-h\left(  t\right)  ,
\end{array}
\right.
\]
the last equation has only the solution $\left(  P^{t,x}\left(  s\right)
,Q^{t,x}\left(  s\right)  \right)  =\left(  -h\left(  t\right)  ,0\right)  ,$
$\forall s\in\left[  t,T\right]  .$

Note that, the corresponding $\mathcal{H}$-function is given by%
\[
\mathcal{H}\left(  s,y,v\right)  =\mathbb{H}\left(  s,y,v,p^{t,x}\left(
s\right)  ,q^{t,x}\left(  s\right)  \right)  =p^{t,x}\left(  s\right)
bv+\sigma q^{t,x}\left(  s\right)  -\frac{1}{2}v^{2},
\]
which is a concave function of $v.$ Then according to the sufficient condition
of optimality, see e.g. Theorem 5.2 pp 138 in \cite{10}, for any fixed initial
pair $\left(  t,x\right)  \in\left[  0,T\right]  \times%
\mathbb{R}
,$ Problem (E) is uniquely solvable with an optimal control $\bar{u}%
^{t,x}\left(  .\right)  $ having the representation%
\[
\bar{u}^{t,x}\left(  s\right)  =b\bar{p}^{t,x}\left(  s\right)  ,\text{
}\forall s\in\left[  t,T\right]  ,
\]
such that the process $\left(  \bar{p}^{t,x}\left(  .\right)  ,\bar{q}%
^{t,x}\left(  .\right)  \right)  $ is the unique adapted solution to the BSDE%
\[
\left\{
\begin{array}
[c]{l}%
d\bar{p}^{t,x}\left(  s\right)  =\bar{q}^{t,x}\left(  s\right)  dW\left(
s\right)  ,\text{ }s\in\left[  t,T\right]  ,\\
\bar{p}^{t,x}\left(  T\right)  =-h\left(  t\right)  \left(  \bar{X}%
^{t,x}\left(  s\right)  -x\right)  .
\end{array}
\right.
\]

By stadard arguments we can show that the processes $\left(  \bar{p}%
^{t,x}\left(  .\right)  ,\bar{q}^{t,x}\left(  .\right)  \right)  $ are
explicitly given by%
\[
\left\{
\begin{array}
[c]{l}%
\bar{p}^{t,x}\left(  s\right)  =-M^{t}\left(  s\right)  \left(  \bar{X}%
^{t,x}\left(  s\right)  -x\right)  ,\text{ }s\in\left[  t,T\right]  ,\\
\bar{q}^{t,x}\left(  s\right)  =-\sigma M^{t}\left(  s\right)  ,\text{ }%
s\in\left[  t,T\right]  ,
\end{array}
\right.
\]
where $\bar{X}^{t,x}\left(  .\right)  $ is the solution of the state equation
corresponding to $\bar{u}^{t,x}\left(  .\right)  ,$ given by%
\[
\left\{
\begin{array}
[c]{l}%
d\bar{X}^{t,x}\left(  s\right)  =b^{2}\bar{p}^{t,x}\left(  s\right)  ds+\sigma
dW\left(  s\right)  ,\text{ }s\in\left[  t,T\right]  ,\\
\bar{X}^{t,x}\left(  t\right)  =x.
\end{array}
\right.
\]
and%
\[
M^{t}\left(  s\right)  =\dfrac{h\left(  t\right)  }{b^{2}h\left(  t\right)
\left(  T-s\right)  +1},\text{ }\forall s\in\left[  t,T\right]  .
\]

A simple computation show that%
\[
\overline{u}^{t,x}\left(  s\right)  =-\dfrac{bh\left(  t\right)  }%
{b^{2}h\left(  t\right)  \left(  T-s\right)  +1}\left(  \bar{X}^{t,x}\left(
s\right)  -x\right)  ,\text{ }\forall s\in\left[  t,T\right]  ,
\]
clearly we have%
\begin{equation}
\overline{u}^{t,x}\left(  s\right)  \neq0,\text{ }\forall s\in\left(
t,T\right]  . \tag{2.5}%
\end{equation}

In the next stage, we will see that Problem (E) is time-inconsistent, for this
we first fix the initial data $\left(  t,x\right)  \in\left[  0,T\right]
\times%
\mathbb{R}
$. Note that, if we assume that the Problem (E) is time-consistent, in the
sense that for any $r\in\left[  t,T\right]  $ the restriction of $\bar
{u}^{t,x}\left(  .\right)  $ on $\left[  r,T\right]  $ is optimal for Problem
(E) with initial pair $\left(  r,\bar{X}^{t,x}\left(  r\right)  \right)  ,$
however as Problem (E) is uniquely solvable for any initial pair, we should
have then $\forall r\in\left(  t,T\right]  $%
\[
\bar{u}^{t,x}\left(  s\right)  =\bar{u}^{r,\bar{X}^{t,x}\left(  r\right)
}\left(  s\right)  =-\dfrac{bh\left(  r\right)  }{b^{2}h\left(  r\right)
\left(  T-s\right)  +1}\left(  \bar{X}^{r,\bar{X}^{t,x}\left(  r\right)
}\left(  s\right)  -\bar{X}^{t,x}\left(  r\right)  \right)  ,\forall
s\in\left[  r,T\right]  ,
\]
where $\bar{X}^{r,\hat{X}^{t,x}\left(  r\right)  }\left(  .\right)  $ solves
the SDE
\[
\left\{
\begin{array}
[c]{l}%
d\bar{X}^{r,\bar{X}^{t,x}\left(  r\right)  }\left(  s\right)  =b^{2}%
\dfrac{h\left(  r\right)  }{b^{2}h\left(  r\right)  \left(  T-s\right)
+1}\left(  \bar{X}^{r,\bar{X}^{t,x}\left(  r\right)  }\left(  s\right)
-\bar{X}^{t,x}\left(  r\right)  \right)  ds+\sigma dW\left(  s\right)
,\forall s\in\left[  r,T\right]  ,\\
\bar{X}^{r,\bar{X}^{t,x}\left(  r\right)  }\left(  r\right)  =\bar{X}%
^{t,x}\left(  r\right)  .
\end{array}
\right.
\]

In particular by the uniqueness of solution to the state SDE we should have%
\[
\bar{u}^{t,x}\left(  r\right)  =-\dfrac{bh\left(  r\right)  }{b^{2}h\left(
r\right)  \left(  T-r\right)  +1}\left(  \bar{X}^{r,\bar{X}^{t,x}\left(
r\right)  }\left(  r\right)  -\bar{X}^{t,x}\left(  r\right)  \right)  =0,
\]
is the only optimal solution of the Problem (E), this contradict $\left(
2.5\right)  $. Therefore, Problem (E) is not time-consistent, and more
precisely, the solution obtained by the classical SMP is wrong and the problem
is rather trivial since the only optimal solution equal to zero.

\section{Characterization of equilibrium strategies}

The purpose of this paper is to characterize open-loop Nash equilibriums
instead of optimal controls. We use the game theoretic approach to handle the
time inconsistency in the same perspective as Ekeland and Lazrak \cite{1},
Bjork and Murgoci \cite{8}. Let us briefly describe the game perspective that
we will consider, as follows.

\begin{itemize}
\item We consider a game with one player at each point $t$ in $\left[
0,T\right]  $. This player represents the incarnation of the controller at
time $t$ and is referred to as \textquotedblleft player $t$\textquotedblright.

\item The $t-th$ player can control the system only at time $t$ by taking
his/her strategy $u\left(  t,.\right)  :\Omega\rightarrow%
\mathbb{R}
^{m}.$

\item A control process $u\left(  .\right)  $ is then viewed as a complete
description of the chosen strategies of all players in the game.

\item The reward to the player $t$ is given by the functional $J\left(
t,\xi,u\left(  .\right)  \right)  $. Note that $J\left(  t,\xi,u\left(
.\right)  \right)  $ depends only on the restriction of the control $u\left(
.\right)  $ to the time interval $\left[  t,T\right]  .$
\end{itemize}

In the above description, we have presented the concept of a \textquotedblleft%
\ Nash equilibrium point\textquotedblright\ of the game. This is an admissible
control process $\hat{u}\left(  .\right)  $ satisfying the following
condition; Suppose that every player $s$, such that $s>t$, will use the
strategy $\hat{u}\left(  s\right)  $. Then the optimal choice for player $t$
is that, he/she also uses the strategy $\hat{u}\left(  t\right)  .$

Nevertheless, the problem with this \textquotedblleft
definition\textquotedblright, is that the individual player $t$ does not
really influence the outcome of the game at all. He/she only chooses the
control at the single point $t,$ and since this is a time set of Lebesgue
measure zero, the control dynamics will not be influenced. Therefore, to
characterize open-loop Nash equilibriums, which have not to be necessary
feedback, we follow \cite{18} who suggest the following formal definition
inspired by \cite{1} and \cite{2}.

Noting that,\textbf{ }for brevity, in the rest of the paper, we suppress the
subscript $\left(  s\right)  $ for the coefficients $A\left(  s\right)
,B\left(  s\right)  ,b\left(  s\right)  ,C_{j}\left(  s\right)  ,D_{j}\left(
s\right)  ,\sigma_{j}\left(  s\right)  $, and we use the notation
$\varrho\left(  z\right)  $ instead of $\varrho\left(  s,z\right)  $\ for
$\varrho=E,F$ and $c.$ In addition, sometimes we simply call $\hat{u}\left(
.\right)  $ an equilibrium control instead of open-loop Nash equilibrium
control when there is no ambiguity.

We define an equilibrium by local spike variation, given an admissible control
$\hat{u}\left(  .\right)  \in\mathcal{L}_{%
\mathcal{F}%
}^{2}\left(  0,T;%
\mathbb{R}
^{m}\right)  .$ For any $t\in\left[  0,T\right]  ,$ $v\in\mathbb{L}^{2}\left(
\Omega,\mathcal{F}_{t},\mathbb{P};%
\mathbb{R}
^{m}\right)  $ and for any $\varepsilon>0,$ define%
\begin{equation}
u^{\varepsilon}\left(  s\right)  =\left\{
\begin{array}
[c]{cc}%
\hat{u}\left(  s\right)  +v, & \text{ for }s\in\left[  t,t+\varepsilon\right)
,\\
\hat{u}\left(  s\right)  , & \text{ for }s\in\left[  t+\varepsilon,T\right]  ,
\end{array}
\right.  \tag{3.1}%
\end{equation}
we have the following definition.

\begin{definition}
[Open-loop Nash equilibrium]An admissible strategy $\hat{u}\left(  .\right)
\in\mathcal{L}_{%
\mathcal{F}%
}^{2}\left(  0,T;%
\mathbb{R}
^{m}\right)  $ is an open-loop Nash equilibrium control for Problem (LQJ) if%
\begin{equation}
\lim_{\varepsilon\downarrow0}\frac{1}{\varepsilon}\left\{  J\left(  t,\hat
{X}\left(  t\right)  ,u^{\varepsilon}\left(  .\right)  \right)  -J\left(
t,\hat{X}\left(  t\right)  ,\hat{u}\left(  .\right)  \right)  \right\}  \geq0,
\tag{3.2}%
\end{equation}
for any $t\in\left[  0,T\right]  ,$ and $v\in\mathbb{L}^{2}\left(
\Omega,\mathcal{F}_{t},\mathbb{P};%
\mathbb{R}
^{m}\right)  .$ The corresponding equilibrium dynamics solves the following
SDE with jumps%
\[
\left\{
\begin{array}
[c]{l}%
d\hat{X}\left(  s\right)  =\left\{  A\hat{X}\left(  s\right)  +B\hat{u}\left(
s\right)  +b\right\}  ds+%
{\displaystyle\sum\limits_{j=1}^{d}}
\left\{  C_{j}\hat{X}\left(  s\right)  +D_{j}\hat{u}\left(  s\right)
+\sigma_{j}\right\}  dW^{j}\left(  s\right) \\
\text{ \ \ \ \ \ \ \ \ \ \ \ \ \ }+%
{\displaystyle\int\limits_{Z}}
\left\{  E\left(  z\right)  \hat{X}\left(  s-\right)  +F\left(  z\right)
\hat{u}\left(  s\right)  +c\left(  z\right)  \right\}  \tilde{N}\left(
ds,dz\right)  ,\text{\ }s\in\left[  0,T\right]  ,\\
\hat{X}_{0}=x_{0}\text{.}%
\end{array}
\right.
\]

\end{definition}

\subsection{The flow of adjoint equations}

First, we introduce the adjoint equations involved in the stochastic maximum
principle which characterize the open-loop Nash equilibrium controls of
Problem (LQJ). Define the Hamiltonian $\mathbb{H}:D\left[  0,T\right]
\times\mathbb{L}^{1}\left(  \Omega,\mathcal{F}_{t},\mathbb{P};%
\mathbb{R}
^{n}\right)  \times%
\mathbb{R}
^{m}\times%
\mathbb{R}
^{n}\times%
\mathbb{R}
^{n\times d}\times\mathbb{L}^{2}\left(  Z,\mathcal{B}\left(  Z\right)  ,\nu;%
\mathbb{R}
^{n}\right)  \rightarrow%
\mathbb{R}
$ by
\begin{equation}%
\begin{array}
[c]{l}%
\mathbb{H}\left(  t,s,X,u,p,q,r\left(  .\right)  \right) \\
=\left\langle p,AX+Bu+b\right\rangle +%
{\displaystyle\sum\limits_{j=1}^{d}}
\left\langle q_{j},D_{j}X+C_{j}u+\sigma_{j}\right\rangle -\dfrac{1}%
{2}\left\langle R\left(  t,s\right)  u,u\right\rangle \\
\text{ \ }+%
{\displaystyle\int\limits_{Z}}
\left\langle r\left(  z\right)  ,E\left(  z\right)  X+F\left(  z\right)
u+c\left(  z\right)  \right\rangle \theta\left(  dz\right)  -\dfrac{1}%
{2}\left(  \left\langle Q\left(  t,s\right)  X,X\right\rangle +\left\langle
\bar{Q}\left(  t,s\right)  \mathbb{E}^{t}\left[  X\right]  ,\mathbb{E}%
^{t}\left[  X\right]  \right\rangle \right)  .
\end{array}
\text{ } \tag{3.3}%
\end{equation}

Let $\hat{u}\left(  .\right)  \in\mathcal{L}_{%
\mathcal{F}%
}^{2}\left(  0,T;%
\mathbb{R}
^{m}\right)  $ and denote by $\hat{X}\left(  .\right)  $ the corresponding
controlled state process. For each $t\in\left[  0,T\right]  $, we introduce
the first order adjoint equation defined on the time interval $\left[
t,T\right]  $, and satisfied by the triple of processes $\left(  p\left(
.;t\right)  ,q\left(  .;t\right)  ,r\left(  .,.;t\right)  \right)  $ as
follows%
\begin{equation}
\left\{
\begin{array}
[c]{l}%
dp\left(  s;t\right)  =-\left\{  A^{\top}p\left(  s;t\right)  +%
{\displaystyle\sum\limits_{j=1}^{d}}
C_{j}^{\top}q_{j}\left(  s;t\right)  \right.  +%
{\displaystyle\int\limits_{Z}}
E\left(  z\right)  ^{\top}r\left(  s,z;t\right)  \theta\left(  dz\right)
-Q\left(  t,s\right)  \hat{X}\left(  s\right) \\
\text{ \ \ \ \ \ \ \ \ \ \ \ \ \ \ \ \ \ \ \ \ \ }\left.  -\bar{Q}\left(
t,s\right)  \mathbb{E}^{t}\left[  \hat{X}\left(  s\right)  \right]  \right\}
ds+%
{\displaystyle\sum\limits_{j=1}^{d}}
q_{j}\left(  s;t\right)  dW^{j}\left(  s\right)  +%
{\displaystyle\int\limits_{Z}}
r\left(  s-,z;t\right)  \tilde{N}\left(  ds,dz\right)  ,\text{\ }s\in\left[
t,T\right]  ,\\
p\left(  T;t\right)  =-G\left(  t\right)  \hat{X}\left(  T\right)  -\bar
{G}\left(  t\right)  \mathbb{E}^{t}\left[  \hat{X}\left(  T\right)  \right]
-\mu_{1}\left(  t\right)  \hat{X}\left(  t\right)  -\mu_{2}\left(  t\right)  ,
\end{array}
\right.  \tag{3.4}%
\end{equation}
where $q\left(  .;t\right)  =\left(  q_{1}\left(  .;t\right)  ,...,q_{d}%
\left(  .;t\right)  \right)  .$ Similarly, we introduce the second order
adjoint equation defined on the time interval $\left[  t,T\right]  ,$ and
satisfied by the triple of processes $\left(  P\left(  .;t\right)
,\Lambda\left(  .;t\right)  ,\Gamma\left(  .,.;t\right)  \right)  $ as follows%
\begin{equation}
\left\{
\begin{array}
[c]{l}%
dP\left(  s;t\right)  =-\left\{  A^{\top}P\left(  s;t\right)  +P\left(
s;t\right)  A+\sum\limits_{j=1}^{d}\left(  C_{j}^{\top}P\left(  s;t\right)
C_{j}\right.  \right.  +\Lambda_{j}\left(  s;t\right)  C_{j}\\
\text{ \ \ \ \ \ \ \ \ \ \ \ \ \ \ \ \ \ \ \ \ \ \ }+\left.  C_{j}^{\top
}\Lambda_{j}\left(  s;t\right)  \right)  +%
{\displaystyle\int_{Z}}
E\left(  z\right)  ^{\top}\left(  \Gamma\left(  s,z;t\right)  +P\left(
s;t\right)  \right)  E\left(  z\right)  \theta\left(  dz\right) \\
\text{ \ \ \ \ \ \ \ \ \ \ \ \ \ \ \ \ \ \ \ \ \ }+%
{\displaystyle\int\limits_{Z}}
\Gamma\left(  s,z;t\right)  E\left(  z\right)  \theta\left(  dz\right)
+\left.
{\displaystyle\int_{Z}}
E\left(  z\right)  ^{\top}\Gamma\left(  s,z;t\right)  \theta\left(  dz\right)
-Q\left(  t,s\right)  \right\}  ds\\
\text{ \ \ \ \ \ \ \ \ \ \ \ \ \ \ \ }+%
{\displaystyle\sum\limits_{j=1}^{d}}
\Lambda_{j}\left(  s;t\right)  dW_{s}^{j}+%
{\displaystyle\int_{Z}}
\Gamma\left(  s-,z;t\right)  \tilde{N}\left(  ds,dz\right)  ,\text{ }%
s\in\left[  t,T\right]  ,\\
P\left(  T;t\right)  =-G\left(  t\right)  ,
\end{array}
\right.  \tag{3.5}%
\end{equation}
where $\Lambda\left(  .;t\right)  =\left(  \Lambda_{1}\left(  .;t\right)
,...,\Lambda_{d}\left(  .;t\right)  \right)  $. Under \textbf{(H1)} the BSDE
$\left(  3.4\right)  $ is uniquely solvable in $\mathcal{S}_{%
\mathcal{F}%
}^{2}\left(  t,T;%
\mathbb{R}
^{n}\right)  \times\mathcal{L}_{%
\mathcal{F}%
}^{2}\left(  t,T;%
\mathbb{R}
^{n\times d}\right)  \times\mathcal{L}_{%
\mathcal{F}%
}^{\theta,2}\left(  \left[  t,T\right]  \times Z;%
\mathbb{R}
^{n}\right)  ,$ see e.g. \cite{23}. Moreover there exists a constant $K>0$
such that%
\begin{equation}
\mathbb{E}\left[  \sup\limits_{t\leq s\leq T}\left\vert p\left(  s;t\right)
\right\vert _{%
\mathbb{R}
^{n}}^{2}\right]  +\mathbb{E}\left[
{\displaystyle\int_{t}^{T}}
\left\vert q\left(  s;t\right)  \right\vert _{%
\mathbb{R}
^{n\times d}}^{2}ds\right]  +\mathbb{E}\left[
{\displaystyle\int_{t}^{T}}
{\displaystyle\int_{Z}}
\left\vert r\left(  s,z;t\right)  \right\vert _{%
\mathbb{R}
^{n}}^{2}\theta\left(  dz\right)  ds\right]  \leq K\left(  1+\left\vert
x_{0}\right\vert ^{2}\right)  . \tag{3.6}%
\end{equation}

In an other hand, noting that the final data of the equation $\left(
3.5\right)  $ is deterministic, it is straightforward to look at a
deterministic solution. \ In addition we have the following representation%
\begin{equation}
\left\{
\begin{array}
[c]{l}%
dP\left(  s;t\right)  =-\left\{  A^{\top}P\left(  s;t\right)  +P\left(
s;t\right)  A+\sum\limits_{j=1}^{d}C_{j}^{\top}P\left(  s;t\right)
C_{j}\right. \\
\text{ \ \ \ \ \ \ \ \ \ \ \ \ \ \ \ \ \ \ \ \ \ \ \ \ }\left.  +%
{\displaystyle\int_{Z}}
E\left(  z\right)  ^{\top}P\left(  s;t\right)  E\left(  z\right)
\theta\left(  dz\right)  -Q\left(  t,s\right)  \right\}  ds,\text{ }%
s\in\left[  t,T\right]  ,\\
P\left(  T;t\right)  =-G\left(  t\right)  ,
\end{array}
\right.  \tag{3.7}%
\end{equation}
which is a uniquely solvable matrix-valued ordinary differential equation.
Indeed, if we define the function $\Phi\left(  s,.\right)  $ for each
$s\in\left[  0,T\right]  ,$ as the fundamental solution of the following
linear SDE%
\begin{equation}
\left\{
\begin{array}
[c]{l}%
d\Phi\left(  s,r\right)  =A\left(  r\right)  \Phi\left(  s,r\right)  dr+%
{\displaystyle\sum\limits_{j=1}^{d}}
C_{j}\left(  r\right)  \Phi\left(  s,r\right)  dW^{j}\left(  r\right)  +%
{\displaystyle\int\limits_{Z}}
E\left(  r,z\right)  \Phi\left(  s,r-\right)  \tilde{N}\left(  dr,dz\right)
,\text{\ }r\in\left[  s,T\right]  ,\\
\Phi\left(  s,s\right)  =I.
\end{array}
\right.  \tag{3.8}%
\end{equation}

Then, by standard arguments based on the Ito's formula we can prove that the
triple $\left(  P\left(  .;t\right)  ,\Lambda\left(  .;t\right)
,\Gamma\left(  .,.;t\right)  \right)  $ solution to $\left(  3.7\right)  $ is
explicitly given by
\begin{equation}
\left\{
\begin{array}
[c]{l}%
P\left(  s;t\right)  =\mathbb{E}^{s}\left[  -\Phi\left(  s,T\right)  ^{\top
}G\left(  t\right)  \Phi\left(  s,T\right)  -%
{\textstyle\int_{s}^{T}}
\Phi\left(  s,r\right)  ^{\top}Q\left(  t,r\right)  \Phi\left(  s,r\right)
ds\right]  ,\text{ }s\in\left[  t,T\right]  ,\\
\Lambda_{j}\left(  s;t\right)  =0,\text{ }s\in\left[  t,T\right]  ,\text{ for
}j=1,2,...,d,\\
\Gamma\left(  s,z;t\right)  =0,\text{ }\left(  s,z\right)  \in\left[
t,T\right]  \times Z.
\end{array}
\right.  \tag{3.9}%
\end{equation}

Next, for each $t\in\left[  0,T\right]  ,$ associated with the 6-tuple
$\left(  \hat{u}\left(  .\right)  ,\hat{X}\left(  .\right)  ,p\left(
.;t\right)  ,q\left(  .,t\right)  ,r\left(  .,.;t\right)  ,P\left(
.;t\right)  \right)  $ we define the $\mathcal{H}_{t}$-function\ as follows%
\begin{align}
\mathcal{H}_{t}\left(  s,X,u\right)   &  =\mathbb{H}\left(  t,s,X,\hat
{u}\left(  s\right)  +u,p\left(  s;t\right)  ,q\left(  s;t\right)  ,r\left(
s,.;t\right)  \right)  +\dfrac{1}{2}%
{\displaystyle\sum\limits_{j=1}^{d}}
u^{\top}D_{j}^{\top}P\left(  s;t\right)  D_{j}u\nonumber\\
&  \text{ \ \ }+\dfrac{1}{2}%
{\displaystyle\int_{Z}}
u^{\top}F\left(  z\right)  ^{\top}P\left(  s;t\right)  F\left(  z\right)
u\theta\left(  dz\right)  , \tag{3.10}%
\end{align}
where $\left(  s,X,u\right)  \in\left[  t,T\right]  \times\mathbb{L}%
^{1}\left(  \Omega,\mathcal{F},\mathbb{P};%
\mathbb{R}
^{n}\right)  \times%
\mathbb{R}
^{m}.$

\subsection{A stochastic maximum principle for equilibrium controls}

In this section, we present a version of Pontryagin's stochastic maximum
principle which characterize the equilibrium controls of Problem (LQJ). We
derive the result by using the second order Taylor expansion in the special
form spike variation $\left(  3.1\right)  $. Here, we don't assume the
non-negativity condition about the matrices $Q,$ $G$ and $R$ as in \cite{18}
and \cite{20}$.$

The following theorem is the first main result of this work, it's providing a
necessary and sufficient condition to characterize the open-loop Nash
equilibrium controls for time-inconsistent Problem (LQJ).

\begin{theorem}
[Stochastic Maximum Principle For Equilibriums]Let \textbf{(H1) }%
holds\textbf{. }Then an admissible control $\hat{u}\left(  .\right)
\in\mathcal{L}_{%
\mathcal{F}%
}^{2}\left(  0,T;%
\mathbb{R}
^{m}\right)  $ is an open-loop Nash equilibrium, if and only if, for any
$t\in\left[  0,T\right]  $, there exist a unique triple of adapted processes
$\left(  p\left(  .;t\right)  ,q\left(  .;t\right)  ,r\left(  .,.;t\right)
\right)  $ which satisfy the BSDE $\left(  3.4\right)  $ and a deterministic
matrix-valued function $P\left(  .;t\right)  $ which satisfies the ODE
$\left(  3.7\right)  $, such that the following condition holds, for all $u\in%
\mathbb{R}
^{m}$%
\begin{equation}%
\begin{array}
[c]{l}%
\mathbb{H}\left(  t,t,\hat{X}\left(  t\right)  ,\hat{u}\left(  t\right)
+u,p\left(  t;t\right)  ,q\left(  t;t\right)  ,r\left(  t,.;t\right)  \right)
-\mathbb{H}\left(  t,t,\hat{X}\left(  t\right)  ,\hat{u}\left(  t\right)
,p\left(  t;t\right)  ,q\left(  t;t\right)  ,r\left(  t,.;t\right)  \right) \\
+\dfrac{1}{2}%
{\displaystyle\sum\limits_{j=1}^{d}}
u^{\top}D_{j}\left(  t\right)  ^{\top}P\left(  t;t\right)  D_{j}\left(
t\right)  u+\dfrac{1}{2}%
{\displaystyle\int\limits_{Z}}
u^{\top}F\left(  t,z\right)  ^{\top}P\left(  t;t\right)  F\left(  t,z\right)
u\theta\left(  dz\right)  \leq0,\text{ }\mathbb{P}-a.s,
\end{array}
\tag{3.11}%
\end{equation}
or equivalently, we have the following two conditions, the first order
equilibrium condition
\begin{equation}
R\left(  t,t\right)  \hat{u}\left(  t\right)  -B\left(  t\right)  ^{\top
}p\left(  t;t\right)  -\sum\limits_{j=1}^{d}D_{j}\left(  t\right)  ^{\top
}q_{j}\left(  t;t\right)  -%
{\displaystyle\int_{Z}}
F\left(  t,z\right)  ^{\top}r\left(  t,z;t\right)  \theta\left(  dz\right)
=0,\text{ }\mathbb{P-}a.s, \tag{3.12}%
\end{equation}
and the second order equilibrium condition%
\begin{equation}
R\left(  t,t\right)  -\sum\limits_{j=1}^{d}D_{j}\left(  t\right)  ^{\top
}P\left(  t;t\right)  D_{j}\left(  t\right)  -%
{\displaystyle\int\limits_{Z}}
F\left(  t,z\right)  ^{\top}P\left(  t;t\right)  F\left(  t,z\right)
\theta\left(  dz\right)  \geq0. \tag{3.13}%
\end{equation}

\end{theorem}

We point that the above result provides a characterisation of open-loop Nash
equilibrium controls via a stochastic maximum principle which is not in the
same setting that the classical stochastic maximum principle for optimal
controls \cite{9} in the sense that, the above result involves the existence
of solutions $\left(  \hat{X}\left(  .\right)  ,\left(  p\left(  .;t\right)
,q\left(  .;t\right)  r\left(  .,.;t\right)  \right)  _{t\in\left[
0,T\right]  }\right)  $ to a "flow" of forward-backward stochastic
differential equations parameterized by $t\in\left[  0,T\right]  $, while the
Pontryagin's stochastic maximum principle for optimal controls involve only
one system of forward-backward stochastic differential equation. Note that for
each $t\in\left[  0,T\right]  $, $\left(  3.4\right)  $ and $\left(
3.5\right)  $ are backward stochastic differential equations. So, as we
consider all $t$ in $\left[  0,T\right]  ,$ all their corresponding adjoint
equations form essentially a "flow" of BSDEs. Moreover, there is an additional
constraint $\left(  3.11\right)  $ which is equivalent to the conditions
$\left(  3.12\right)  $ and $\left(  3.13\right)  $ that acts on the flow only
when $s=t$.

Our goal now, is to give a proof of the Theorem 3.2. The main idea is still
based on the variational techniques in the same spirit of proving the
stochastic Pontryagin's maximum principle \cite{15}.

Let $\hat{u}\left(  .\right)  \in\mathcal{L}_{%
\mathcal{F}%
}^{2}\left(  0,T;%
\mathbb{R}
^{m}\right)  $ be an admissible control and $\hat{X}\left(  .\right)  $ the
corresponding controlled process solution to the state equation$.$ Consider
the perturbed control $u^{\varepsilon}\left(  .\right)  $ defined by the spike
variation $\left(  3.1\right)  $ for some fixed arbitrary $t\in\left[
0,T\right]  ,$ $v\in\mathbb{L}^{2}\left(  \Omega,\mathcal{F}_{t},\mathbb{P};%
\mathbb{R}
^{m}\right)  $ and $\varepsilon\in\left[  0,T-t\right]  .$ Denote by $\hat
{X}^{\varepsilon}\left(  .\right)  $ the solution of the state equation
corresponding to $u^{\varepsilon}\left(  .\right)  $. Since the coefficients
of the controlled state equation are linear, then by the standard perturbation
approach, see e.g. \cite{15}, we have%
\begin{equation}
\hat{X}^{\varepsilon}\left(  s\right)  -\hat{X}\left(  s\right)
=y^{\varepsilon,v}\left(  s\right)  +z^{\varepsilon,v}\left(  s\right)
,\text{ }s\in\left[  t,T\right]  , \tag{3.14}%
\end{equation}
where $y^{\varepsilon,v}\left(  .\right)  $ and $z^{\varepsilon,v}\left(
.\right)  $ solve the following linear stochastic differential equations,
respectively%
\begin{equation}
\left\{
\begin{array}
[c]{l}%
dy^{\varepsilon,v}\left(  s\right)  =Ay^{\varepsilon,v}\left(  s\right)
ds+\sum\limits_{j=1}^{d}\left\{  C_{j}y^{\varepsilon,v}\left(  s\right)
+D_{j}v1_{\left[  t,t+\varepsilon\right)  }\left(  s\right)  \right\}
dW^{j}\left(  s\right) \\
\text{ \ \ \ \ \ \ \ \ \ \ \ \ \ }+%
{\displaystyle\int_{Z}}
\left\{  E\left(  z\right)  y^{\varepsilon,v}\left(  s-\right)  +F\left(
z\right)  v1_{\left[  t,t+\varepsilon\right)  }\left(  s\right)  \right\}
\tilde{N}\left(  ds,dz\right)  ,\text{ }s\in\left[  t,T\right]  ,\\
y^{\varepsilon,v}\left(  t\right)  =0,\text{\ }%
\end{array}
\right.  \tag{3.15}%
\end{equation}
and%
\begin{equation}
\left\{
\begin{array}
[c]{l}%
dz^{\varepsilon,v}\left(  s\right)  =\left\{  Az^{\varepsilon,v}\left(
s\right)  +Bv1_{\left[  t,t+\varepsilon\right)  }\left(  s\right)  \right\}
ds+\sum\limits_{j=1}^{d}C_{j}z^{\varepsilon,v}\left(  s\right)  dW^{j}\left(
s\right) \\
\text{ \ \ \ \ \ \ \ \ \ \ \ \ \ \ \ }+%
{\displaystyle\int_{Z}}
E\left(  z\right)  z^{\varepsilon,v}\left(  s\right)  \tilde{N}\left(
ds,dz\right)  ,\text{\ }s\in\left[  t,T\right]  ,\\
z^{\varepsilon,v}\left(  t\right)  =0.
\end{array}
\right.  \tag{3.16}%
\end{equation}

Now, we present the following technical lemma needed later in this study, see
the Appendix \textbf{A.1.} for its proof.

\begin{lemma}
Under assumption \textbf{(H1)},\ the following estimates hold%
\begin{align}
&  \mathbb{E}^{t}\left[  y^{\varepsilon}\left(  s\right)  \right]  =0,\text{
a.e. }s\in\left[  t,T\right]  \text{ and }\sup\limits_{s\in\left[  t,T\right]
}\left\vert \mathbb{E}^{t}\left[  z^{\varepsilon}\left(  s\right)  \right]
\right\vert ^{2}=O\left(  \varepsilon^{2}\right)  ,\tag{3.17}\\
&  \mathbb{E}^{t}\left[  \sup\limits_{s\in\left[  t,T\right]  }\left\vert
y^{\varepsilon}\left(  s\right)  \right\vert ^{2}\right]  =O\left(
\varepsilon\right)  \text{ and }\mathbb{E}^{t}\left[  \sup\limits_{s\in\left[
t,T\right]  }\left\vert z^{\varepsilon}\left(  s\right)  \right\vert
^{2}\right]  =O\left(  \varepsilon^{2}\right)  . \tag{3.18}%
\end{align}

Moreover, we have the equality%
\begin{equation}%
\begin{array}
[c]{l}%
J\left(  t,\hat{X}\left(  t\right)  ,u^{\varepsilon}\left(  .\right)  \right)
-J\left(  t,\hat{X}\left(  t\right)  ,\hat{u}\left(  .\right)  \right) \\
=-\mathbb{E}^{t}\left[
{\displaystyle\int\limits_{t}^{T}}
\left\{  \mathbb{H}\left(  t,s,\hat{X}\left(  s\right)  ,\hat{u}\left(
s\right)  +v,p\left(  s;t\right)  ,q\left(  s;t\right)  ,r\left(
s,.;t\right)  \right)  \right.  \right. \\
\text{ \ \ \ \ }\ \ \ \ \ \ \ \ \text{\ }-\left.  \mathbb{H}\left(
t,s,\hat{X}\left(  s\right)  ,\hat{u}\left(  s\right)  ,p\left(  s;t\right)
,q\left(  s;t\right)  ,r\left(  s,.;t\right)  \right)  \right\}  1_{\left[
t,t+\varepsilon\right)  }\left(  s\right)  ds\\
\text{ \ \ \ \ \ \ \ \ \ }+\dfrac{1}{2}%
{\displaystyle\int\limits_{t}^{T}}
\left\{  \sum\limits_{j=1}^{d}v^{\top}D_{j}^{\top}P\left(  s;t\right)
D_{j}v\right.  +\left.  \left.
{\displaystyle\int_{Z}}
v^{\top}F\left(  z\right)  ^{\top}P\left(  s;t\right)  F\left(  z\right)
v\theta\left(  dz\right)  \right\}  1_{\left[  t,t+\varepsilon\right)
}\left(  s\right)  ds\right]  +o\left(  \varepsilon\right)  .
\end{array}
\tag{3.19}%
\end{equation}

\end{lemma}

Now, we are ready to give the proof of the Theorem 3.2.

\textbf{Proof of Theorem 3.2.} Given an open-loop Nash equilibrium $\hat
{u}\left(  .\right)  $, then for any $t\in\left[  0,T\right]  $ and
$v\in\mathbb{L}^{2}\left(  \Omega,\mathcal{F}_{t},\mathbb{P};%
\mathbb{R}
^{m}\right)  ,$ we have clearly%
\[
\lim_{\varepsilon\downarrow0}\frac{1}{\varepsilon}\left\{  J\left(  t,\hat
{X}\left(  t\right)  ,\hat{u}\left(  .\right)  \right)  -J\left(  t,\hat
{X}\left(  t\right)  ,u^{\varepsilon}\left(  .\right)  \right)  \right\}
\leq0,
\]
from which we deduce

$\underset{\varepsilon\rightarrow0}{\lim}\dfrac{1}{\varepsilon}\mathbb{E}%
^{t}\left[
{\displaystyle\int_{t}^{T}}
\left\{  \mathbb{H}\left(  t,s,\hat{X}\left(  s\right)  ,\hat{u}\left(
s\right)  +v,p\left(  s;t\right)  ;q\left(  s;t\right)  ,r\left(
s,.;t\right)  \right)  \right.  \right.  $

$\ \ \ \ \ \ \ \ \ \ \ \ \ \ \ \ \ \ \ \ \ \ \ -\mathbb{H}\left(  t,s,\hat
{X}\left(  s\right)  ,\hat{u}\left(  s\right)  ,p\left(  s;t\right)  ,q\left(
s;t\right)  ,r\left(  s,.;t\right)  \right)  $

$\ \ \ \ \ \ \ \ \ \ \ \ \ \ \ \ \ \ \ \ \ +\dfrac{1}{2}%
{\displaystyle\sum\limits_{j=1}^{d}}
v^{\top}D_{j}^{\top}P\left(  s;t\right)  D_{j}v+\left.  \left.  \dfrac{1}{2}%
{\displaystyle\int_{Z}}
v^{\top}F\left(  z\right)  ^{\top}P\left(  s;t\right)  F\left(  z\right)
v\theta\left(  dz\right)  \right\}  1_{\left[  t,t+\varepsilon\right)
}\left(  s\right)  ds\right]  \leq0,$

which leads to%
\[%
\begin{array}
[c]{c}%
\mathbb{H}\left(  t,t,\hat{X}\left(  t\right)  ,\hat{u}\left(  t\right)
+v,p\left(  t;t\right)  ,q\left(  t;t\right)  ,r\left(  t,.;t\right)  \right)
-\mathbb{H}\left(  t,t,\hat{X}\left(  t\right)  ,\hat{u}\left(  t\right)
,p\left(  t;t\right)  ,q\left(  t;t\right)  ,r\left(  t,.;t\right)  \right) \\
\text{ \ \ \ \ \ \ }+\dfrac{1}{2}%
{\displaystyle\sum\limits_{j=1}^{d}}
v^{\top}D_{j}\left(  t\right)  ^{\top}P\left(  t;t\right)  D_{j}\left(
t\right)  v+\dfrac{1}{2}%
{\displaystyle\int_{Z}}
v^{\top}F\left(  t,z\right)  ^{\top}P\left(  t;t\right)  F\left(  t,z\right)
v\theta\left(  dz\right)  \leq0,\text{ }\mathbb{P}-a.s,
\end{array}
\]

Therefore, by stetting $v\equiv u$ for an arbitrarily $u\in%
\mathbb{R}
^{m}$ we obtain $\left(  3.11\right)  .$

Conversely, given an admissible control $\hat{u}\left(  .\right)
\in\mathcal{L}_{%
\mathcal{F}%
}^{2}\left(  0,T;%
\mathbb{R}
^{m}\right)  .$ Suppose that for any $t\in\left[  0,T\right]  ,$ the
variational inequality $\left(  3.11\right)  $ holds. Then for any
$v\in\mathbb{L}^{2}\left(  \Omega,\mathcal{F}\left(  t\right)  ,\mathbb{P};%
\mathbb{R}
^{m}\right)  $ it yields%
\[%
\begin{array}
[c]{c}%
\mathbb{H}\left(  t,t,\hat{X}\left(  t\right)  ,\hat{u}\left(  t\right)
+v,p\left(  t;t\right)  ,q\left(  t;t\right)  ,r\left(  t,.;t\right)  \right)
-\mathbb{H}\left(  t,t,\hat{X}\left(  t\right)  ,\hat{u}\left(  t\right)
,p\left(  t;t\right)  ,q\left(  t;t\right)  ,r\left(  t,.;t\right)  \right) \\
+\dfrac{1}{2}%
{\displaystyle\sum\limits_{j=1}^{d}}
v^{\top}D_{j}\left(  t\right)  ^{\top}P\left(  t;t\right)  D_{j}\left(
t\right)  v+\dfrac{1}{2}%
{\displaystyle\int_{Z}}
v^{\top}F\left(  t,z\right)  ^{\top}P\left(  t;t\right)  F\left(  t,z\right)
v\theta\left(  dz\right)  \leq0,\text{ }\mathbb{P}-a.s,
\end{array}
\]
consequently

$\lim\limits_{\varepsilon\downarrow0}\dfrac{1}{\varepsilon}\mathbb{E}%
^{t}\left[
{\displaystyle\int_{t}^{t+\varepsilon}}
\left\{  \mathbb{H}\left(  t,s,\hat{X}\left(  s\right)  ,\hat{u}\left(
s\right)  +v,p\left(  s;t\right)  ,q\left(  s;t\right)  ,r\left(
s,.;t\right)  \right)  \right.  \right.  $

$\ \ \ \ \ \ \ \ \ \ \ \ \ \ \ \ \ \ \ \ \ \ \ \ -\mathbb{H}\left(
t,s,\hat{X}\left(  s\right)  ,\hat{u}\left(  s\right)  ,p\left(  s;t\right)
,q\left(  s;t\right)  ,r\left(  s,.;t\right)  \right)  $

$\ \ \ \ \ \ \ \ \ \ \ \ \ \ \ \ \ \ \ \ \ $\ $+\dfrac{1}{2}%
{\displaystyle\sum\limits_{j=1}^{d}}
v^{\top}D_{j}^{\top}P\left(  s;t\right)  D_{j}v+\left.  \left.  \dfrac{1}{2}%
{\displaystyle\int_{Z}}
v^{\top}F\left(  z\right)  ^{\top}P\left(  s;t\right)  F\left(  z\right)
v\theta\left(  dz\right)  \right\}  ds\right]  \leq0.$

Hence%
\[
\lim_{\varepsilon\downarrow0}\frac{1}{\varepsilon}\left\{  J\left(  t,\hat
{X}\left(  t\right)  ,\hat{u}\left(  .\right)  \right)  -J\left(  t,\hat
{X}\left(  t\right)  ,u^{\varepsilon}\left(  .\right)  \right)  \right\}
\leq0,
\]
thus $\hat{u}\left(  .\right)  $ is an equilibrium control.

Easy manipulations show that the variational inequality $\left(  3.11\right)
$ is equivalent to%
\[
\mathcal{H}_{t}\left(  t,\hat{X}\left(  t\right)  ,0\right)  =\max_{u\in%
\mathbb{R}
^{m}}\mathcal{H}_{t}\left(  t,\hat{X}\left(  t\right)  ,u\right)  ,
\]
then $\left(  3.12\right)  $ and $\left(  3.13\right)  $ follow respectively
from the following first order and second order conditions of the maximum
point $\hat{u}=0$ for the quadratic function $u\rightarrow\mathcal{H}%
_{t}\left(  t,\hat{X}\left(  t\right)  ,.\right)  $
\[
\dfrac{\partial\mathcal{H}_{t}}{\partial u}\left(  t,\hat{X}\left(  t\right)
,0\right)  =0\text{ and }\dfrac{\partial^{2}\mathcal{H}_{t}}{\partial u^{2}%
}\left(  t,\hat{X}\left(  t\right)  ,u\right)  \leq0.
\]

Then, the required result directly follows.\eop

In Theorem 3.2, in view of condition $\left(  3.11\right)  $, as long as the
term%
\[
-\left(  \sum\limits_{j=1}^{d}D_{j}\left(  t\right)  ^{\top}P\left(
t;t\right)  D_{j}\left(  t\right)  +\right.  \left.
{\displaystyle\int\limits_{Z}}
F\left(  t,z\right)  ^{\top}P\left(  t;t\right)  F\left(  t,z\right)
\theta\left(  dz\right)  \right)  ,
\]
for each $t\in\left[  0,T\right]  $ is sufficiently positive definite, the
necessary and sufficient condition for equilibriums might still be satisfied
even if $R\left(  t,t\right)  $ is negative. This is different from \cite{18}
and \cite{20} where the authors have assumed the non-negativity of the
matrices $Q,$ $G$ and $R$ in order to state their stochastic maximum principle
for open-loop Nash equilibriums. Moreover, in view of $\left(  3.9\right)  $
in the case where $Q\left(  t,s\right)  \geq0$ for every $s\in\left[
t,T\right]  ,$ and $G\left(  t\right)  \geq0,$ it follows that the solution of
the second order adjoint equation satisfies $P\left(  t;t\right)  \leq0,$ then
if further we have $R\left(  t,t\right)  \geq0,$ Thus the condition that
\[
R\left(  t,t\right)  -\sum\limits_{j=1}^{d}D_{j}\left(  t\right)  ^{\top
}P\left(  t;t\right)  D_{j}\left(  t\right)  -%
{\displaystyle\int\limits_{Z}}
F\left(  t,z\right)  ^{\top}P\left(  t;t\right)  F\left(  t,z\right)
\theta\left(  dz\right)  \geq0,
\]
is obviously satisfied. Therefore, we summarize the main theorem into the
following Corollary.

\begin{corollary}
Let \textbf{(H1)-(H2) }hold\textbf{. }Then an admissible control $\hat
{u}\left(  .\right)  \in\mathcal{L}_{%
\mathcal{F}%
}^{2}\left(  0,T;%
\mathbb{R}
^{m}\right)  $ is an equilibrium control, if and only if, for any $t\in\left[
0,T\right]  $, there exists a triple of adapted processes $\left(  p\left(
.;t\right)  ,q\left(  .;t\right)  ,r\left(  .,.;t\right)  \right)  $ which
satisfies the BSDE $\left(  3.4\right)  ,$ with only the first order condition
$\left(  3.12\right)  $ holds.
\end{corollary}

\section{Linear feedback stochastic equilibrium control}

In this section, we consider only the case where the Brownian motion is
one-dimensional $(d=1)$ for simplicity of presentation. There is no essential
difficulty with the multidimensional Brownian motions. All the indices $j$
will then be dropped. Our goal is to obtain a state feedback representation of
an equilibrium control for Problem (LQJ) via some class of ordinary
differential equations. Suppose that $\hat{u}\left(  .\right)  $ is an
equilibrium control and denote by $\hat{X}\left(  .\right)  $ the
corresponding controlled process$.$ Then in view of Theorem $3.2,$ there
exists a flow of triple of adapted processes $\left(  p\left(  .;t\right)
,q\left(  .;t\right)  ,r\left(  .,.;t\right)  \right)  _{t\in\left[
0,T\right)  }$ for which the 3-tuple $\left(  \hat{u}\left(  .\right)
,\hat{X}\left(  .\right)  ,\left(  p\left(  .;t\right)  ,q\left(  .;t\right)
,r\left(  .,.;t\right)  \right)  _{t\in\left[  0,T\right)  }\right)  $ solves
the following flow of forward-backward SDE with jumps, parametrized by
$t\in\left[  0,T\right]  $%
\begin{equation}
\left\{
\begin{array}
[c]{l}%
d\hat{X}\left(  s\right)  =\left\{  A\hat{X}\left(  s\right)  +B\hat{u}\left(
s\right)  +b\right\}  ds+\left\{  C\hat{X}\left(  s\right)  +D\hat{u}\left(
s\right)  +\sigma\right\}  dW\left(  s\right) \\
\text{ \ \ \ \ \ \ \ \ \ \ \ \ \ }+%
{\displaystyle\int\limits_{Z}}
\left\{  E\left(  z\right)  \hat{X}\left(  s-\right)  +F\left(  z\right)
\hat{u}\left(  s\right)  +c\left(  z\right)  \right\}  \tilde{N}\left(
ds,dz\right)  ,\text{\ }s\in\left[  0,T\right]  ,\\
dp\left(  s;t\right)  =-\left\{  A^{\top}p\left(  s;t\right)  +C^{\top
}q\left(  s;t\right)  \right.  +%
{\displaystyle\int\limits_{Z}}
E\left(  z\right)  ^{\top}r\left(  s,z;t\right)  \theta\left(  dz\right)
-Q\left(  t,s\right)  \hat{X}\left(  s\right) \\
\text{\ \ \ \ \ \ \ \ \ \ \ \ \ \ \ \ \ \ \ \ }\left.  -\bar{Q}\left(
t,s\right)  \mathbb{E}^{t}\left[  \hat{X}\left(  s\right)  \right]  \right\}
ds+q\left(  s;t\right)  dW\left(  s\right)  +%
{\displaystyle\int\limits_{Z}}
r\left(  s-,z;t\right)  \tilde{N}\left(  ds,dz\right)  ,\ 0\leq t\leq s\leq
T,\\
\hat{X}_{0}=x_{0}\text{, }p\left(  T;t\right)  =-G\left(  t\right)  \hat
{X}\left(  T\right)  -\bar{G}\left(  t\right)  \mathbb{E}^{t}\left[  \hat
{X}\left(  T\right)  \right]  -\mu_{1}\left(  t\right)  \hat{X}\left(
t\right)  -\mu_{2}\left(  t\right)  ,\text{ }t\in\left[  0,T\right]  ,
\end{array}
\right.  \tag{4.1}%
\end{equation}
with the condition
\begin{equation}
R\left(  t,t\right)  \hat{u}\left(  t\right)  -B^{\top}p\left(  t;t\right)
-D\left(  t\right)  ^{\top}q\left(  t;t\right)  -%
{\displaystyle\int_{Z}}
F\left(  t,z\right)  ^{\top}r\left(  t,z;t\right)  \theta\left(  dz\right)
=0,\text{ }\mathbb{P}-a.s,\text{ }\forall t\in\left[  0,T\right]  . \tag{4.2}%
\end{equation}

Now, to solve the above stochastic system$,$ we conjecture that $\hat
{X}\left(  .\right)  $ and $p\left(  .;t\right)  $ for $t\in\left[
0,T\right)  $\ are related by the following relation%
\begin{equation}
p\left(  s;t\right)  =-M\left(  t,s\right)  \hat{X}\left(  s\right)  -\bar
{M}\left(  t,s\right)  \mathbb{E}^{t}\left[  \hat{X}\left(  s\right)  \right]
-\Upsilon\left(  t,s\right)  \hat{X}\left(  t\right)  -\varphi\left(
t,s\right)  ,\text{ }\left(  t,s\right)  \in D\left[  0,T\right]  , \tag{4.3}%
\end{equation}
for some deterministic functions $M\left(  .,.\right)  ,\bar{M}\left(
.,.\right)  ,\Upsilon\left(  .,.\right)  \in C^{0,1}\left(  D\left[
0,T\right]  ,%
\mathbb{R}
^{n\times n}\right)  $ and $\varphi\left(  .,.\right)  \in C^{0,1}\left(
D\left[  0,T\right]  ,%
\mathbb{R}
^{n}\right)  $ such that
\begin{equation}
M\left(  t,T\right)  =G\left(  t\right)  ,\text{ }\bar{M}\left(  t,T\right)
=\bar{G}\left(  t\right)  ,\text{ }\Upsilon\left(  t,T\right)  =\mu_{1}\left(
t\right)  \text{, }\varphi\left(  t,T\right)  =\mu_{2}\left(  t\right)
,\text{ }t\in\left[  0,T\right]  . \tag{4.4}%
\end{equation}

Applying Ito's formula to $\left(  4.3\right)  $ and using $\left(
4.1\right)  $, it yields%
\begin{align}
dp\left(  s;t\right)   &  =\left\{  -\dfrac{\partial M}{\partial s}\left(
t,s\right)  \hat{X}\left(  s\right)  -\dfrac{\partial\bar{M}}{\partial
s}\left(  t,s\right)  \mathbb{E}^{t}\left[  \hat{X}\left(  s\right)  \right]
-\dfrac{\partial\Upsilon}{\partial s}\left(  t,s\right)  \hat{X}\left(
t\right)  -\dfrac{\partial\varphi}{\partial s}\left(  t,s\right)  \right.
\nonumber\\
&  \text{ \ \ \ \ \ \ }-\left.  M\left(  t,s\right)  \left(  A\hat{X}\left(
s\right)  +Bu\left(  s\right)  +b\right)  -\bar{M}\left(  t,s\right)  \left(
A\mathbb{E}^{t}\left[  \hat{X}\left(  s\right)  \right]  +B\mathbb{E}%
^{t}\left[  u\left(  s\right)  \right]  +b\right)  \right\}  ds\nonumber\\
&  \text{ \ \ }-M\left(  t,s\right)  \left(  C\hat{X}\left(  s\right)
+D\hat{u}\left(  s\right)  +\sigma\right)  dW\left(  s\right) \nonumber\\
&  \text{ \ \ }-\left.  \int_{Z}M\left(  t,s\right)  \left(  E\left(
z\right)  \hat{X}\left(  s-\right)  +F\left(  z\right)  \hat{u}\left(
s\right)  +c\left(  z\right)  \right)  \tilde{N}\left(  ds,dz\right)  \right.
,\nonumber\\
&  =-\left\{  A^{\top}p\left(  s;t\right)  +C^{\top}q\left(  s;t\right)
+\int_{Z}E\left(  z\right)  ^{\top}r\left(  s,z;t\right)  \theta\left(
dz\right)  \right.  -Q\left(  t,s\right)  \hat{X}\left(  s\right) \nonumber\\
&  \text{ \ \ \ \ \ \ \ \ \ \ }\left.  -\bar{Q}\left(  t,s\right)
\mathbb{E}^{t}\left[  \hat{X}\left(  s\right)  \right]  \right\}  ds+q\left(
s;t\right)  dW\left(  s\right)  +\int_{Z}r\left(  s-,z;t\right)  \tilde
{N}\left(  ds,dz\right)  ,\text{ }s\in\left[  t,T\right]  , \tag{4.5}%
\end{align}
and we obtain
\begin{align}
&  q\left(  s;t\right)  =-M\left(  t,s\right)  \left(  C\hat{X}\left(
s\right)  +D\hat{u}\left(  s\right)  +\sigma\right)  ,\tag{4.6}\\
&  r\left(  s,z;t\right)  =-M\left(  t,s\right)  \left(  E\left(  z\right)
\hat{X}\left(  s\right)  +F\left(  z\right)  \hat{u}\left(  s\right)
+c\left(  z\right)  \right)  . \tag{4.7}%
\end{align}

By taking $\left(  4.6\right)  $ and $\left(  4.7\right)  $ in $\left(
4.2\right)  $ we get
\begin{align*}
0  &  =R\left(  t,t\right)  \hat{u}\left(  t\right)  +B\left(  t\right)
^{\top}\left(  \left(  M\left(  t,t\right)  +\bar{M}\left(  t,t\right)
+\Upsilon\left(  t,t\right)  \right)  \hat{X}\left(  t\right)  +\varphi\left(
t,t\right)  \right) \\
&  \text{ \ \ \ }+D\left(  t\right)  ^{\top}M\left(  t,t\right)  \left(
C\left(  t\right)  \hat{X}\left(  t\right)  +D\left(  t\right)  \hat{u}\left(
t\right)  +\sigma\left(  t\right)  \right) \\
&  \text{ \ \ \ }+\int_{Z}F\left(  t,z\right)  ^{\top}M\left(  t,t\right)
\left(  E\left(  t,z\right)  \hat{X}\left(  t\right)  +F\left(  t,z\right)
\hat{u}\left(  t\right)  +c\left(  t,z\right)  \right)  \theta\left(
dz\right)  ,\\
&  =\left\{  R\left(  t,t\right)  +D\left(  t\right)  ^{\top}M\left(
t,t\right)  D\left(  t\right)  +\int_{Z}F\left(  t,z\right)  ^{\top}M\left(
t,t\right)  F\left(  t,z\right)  \theta\left(  dz\right)  \right\}  \hat
{u}\left(  t\right) \\
&  \text{ \ \ \ \ }+\left\{  B\left(  t\right)  ^{\top}\left(  M\left(
t,t\right)  +\bar{M}\left(  t,t\right)  +\Upsilon\left(  t,t\right)  \right)
+D\left(  t\right)  ^{\top}M\left(  t,t\right)  C\left(  t\right)  \right. \\
&  \text{ \ \ \ \ \ \ \ \ }\left.  +\int_{Z}F\left(  t,z\right)  ^{\top
}M\left(  t,t\right)  E\left(  t,z\right)  \theta\left(  dz\right)  \right\}
\hat{X}\left(  t\right) \\
&  \text{ \ \ \ \ }+B\left(  t\right)  ^{\top}\varphi\left(  t,t\right)
+D\left(  t\right)  ^{\top}M\left(  t,t\right)  \sigma\left(  t\right)
+\int_{Z}F\left(  t,z\right)  ^{\top}M\left(  t,t\right)  c\left(  t,z\right)
\theta\left(  dz\right)  .
\end{align*}

Thus if we assume that $\Theta\left(  t\right)  =\left(  R\left(  t,t\right)
+D\left(  t\right)  ^{\top}M\left(  t,t\right)  D\left(  t\right)  +\int
_{Z}F\left(  t,z\right)  ^{\top}M\left(  t,t\right)  F\left(  t,z\right)
\nu\left(  dz\right)  \right)  ^{-1}$ exists, then we deduce that $\hat
{u}\left(  .\right)  $ admits the following feedback representation%
\begin{equation}
\hat{u}\left(  t\right)  =-\Psi\left(  t\right)  \hat{X}\left(  t\right)
-\psi\left(  t\right)  , \tag{4.8}%
\end{equation}
where $\Psi\left(  t\right)  $ and $\psi\left(  t\right)  $ are given by
\begin{equation}
\left\{
\begin{array}
[c]{l}%
\Psi\left(  t\right)  =\Theta\left(  t\right)  \left(  B\left(  t\right)
^{\top}\left(  M\left(  t,t\right)  +\bar{M}\left(  t,t\right)  +\Upsilon
\left(  t,t\right)  \right)  \right.  +D\left(  t\right)  ^{\top}M\left(
t,t\right) \\
\text{ \ \ \ \ \ \ \ \ \ \ \ \ \ \ \ \ \ \ \ \ \ \ }+%
{\displaystyle\int_{Z}}
\left.  F\left(  t,z\right)  ^{\top}M\left(  t,t\right)  E\left(  t,z\right)
\theta\left(  dz\right)  \right)  ,\\
\psi\left(  t\right)  =\Theta\left(  t\right)  \left(  B\left(  t\right)
^{\top}\varphi\left(  t,t\right)  +D\left(  t\right)  ^{\top}M\left(
t,t\right)  \sigma\left(  t\right)  +%
{\displaystyle\int_{Z}}
F\left(  t,z\right)  ^{\top}M\left(  t,t\right)  c\left(  t,z\right)
\theta\left(  dz\right)  \right)  .
\end{array}
\right.  \tag{4.9}%
\end{equation}

Therefore, for any $\left(  t,s\right)  \in D\left[  0,T\right]  ,$ we have
\begin{equation}
\mathbb{E}^{t}\left[  \hat{u}\left(  s\right)  \right]  =\Psi\left(  s\right)
\mathbb{E}^{t}\left[  \hat{X}\left(  s\right)  \right]  +\psi\left(  s\right)
. \tag{4.10}%
\end{equation}

Next, comparing the $ds$ term in $\left(  4.5\right)  $ by the one in $\left(
4.1\right)  ,$ then by using the expressions $\left(  4.8\right)  $ and
$\left(  4.10\right)  $, we obtain%
\begin{align*}
0  &  =\left\{  \dfrac{\partial M}{\partial s}\left(  t,s\right)  +M\left(
t,s\right)  A+A^{\top}M\left(  t,s\right)  +C^{\top}M\left(  t,s\right)
C\right.  +%
{\displaystyle\int_{Z}}
E\left(  z\right)  ^{\top}M\left(  t,s\right)  E\left(  z\right)
\theta\left(  dz\right) \\
&  \text{ \ \ \ \ \ \ }-\left.  \left(  M\left(  t,s\right)  B+C^{\top
}M\left(  t,s\right)  D+%
{\displaystyle\int_{Z}}
E\left(  z\right)  ^{\top}M\left(  t,s\right)  F\left(  z\right)
\theta\left(  dz\right)  \right)  \Psi\left(  s\right)  +Q\left(  t,s\right)
\right\}  \hat{X}\left(  s\right) \\
&  \text{ \ }+\left\{  \dfrac{\partial\bar{M}}{\partial s}\left(  t,s\right)
+\bar{M}\left(  t,s\right)  A+A^{\top}\bar{M}\left(  t,s\right)  -\bar
{M}\left(  t,s\right)  B\Psi\left(  s\right)  +\bar{Q}\left(  t,s\right)
\right\}  \mathbb{E}^{t}\left[  \hat{X}\left(  s\right)  \right] \\
&  \text{ \ }+\left\{  \dfrac{\partial\Upsilon}{\partial s}\left(  t,s\right)
+A^{\top}\Upsilon\left(  t,s\right)  \right\}  \hat{X}\left(  t\right) \\
&  \text{ \ }+\dfrac{\partial\varphi}{\partial s}\left(  t,s\right)  +\left(
M\left(  t,s\right)  +\bar{M}\left(  t,s\right)  \right)  \left(
b-B\psi\left(  s\right)  \right)  +A^{\top}\varphi\left(  t,s\right) \\
&  \text{ \ }+C^{\top}M\left(  t,s\right)  \left(  \sigma-D\psi\left(
s\right)  \right)  +%
{\displaystyle\int_{Z}}
E\left(  z\right)  ^{\top}M\left(  t,s\right)  \left(  c\left(  z\right)
-F\left(  z\right)  \psi\left(  s\right)  \right)  \theta\left(  dz\right)
\end{align*}

This suggests that the functions $M\left(  .,.\right)  ,\bar{M}\left(
.,.\right)  ,\Upsilon\left(  .,.\right)  $ and $\varphi\left(  .,.\right)  $
solve the following system of ordinary differential equations, for $\left(
t,s\right)  \in D\left[  0,T\right]  $%
\begin{equation}
\left\{
\begin{array}
[c]{l}%
0=\dfrac{\partial M}{\partial s}\left(  t,s\right)  +M\left(  t,s\right)
A+A^{\top}M\left(  t,s\right)  +C^{\top}M\left(  t,s\right)  C+%
{\displaystyle\int_{Z}}
E\left(  z\right)  ^{\top}M\left(  t,s\right)  E\left(  z\right)
\theta\left(  dz\right) \\
\text{ \ \ \ \ \ \ }-\left(  M\left(  t,s\right)  B+C^{\top}M\left(
t,s\right)  D+%
{\displaystyle\int_{Z}}
E\left(  z\right)  ^{\top}M\left(  t,s\right)  F\left(  z\right)
\theta\left(  dz\right)  \right)  \Psi\left(  s\right)  +Q\left(  t,s\right)
,\\
0=\left.  \dfrac{\partial\bar{M}}{\partial s}\left(  t,s\right)  +\bar
{M}\left(  t,s\right)  A+A^{\top}\bar{M}\left(  t,s\right)  -\bar{M}\left(
t,s\right)  B\Psi\left(  s\right)  +\bar{Q}\left(  t,s\right)  ,\right.
\text{ }\\
0=\left.  \dfrac{\partial\Upsilon}{\partial s}\left(  t,s\right)  +A^{\top
}\Upsilon\left(  t,s\right)  \right.  ,\\
0=\dfrac{\partial\varphi}{\partial s}\left(  t,s\right)  +\left(  M\left(
t,s\right)  +\bar{M}\left(  t,s\right)  \right)  \left(  b-B\psi\left(
s\right)  \right)  +A^{\top}\varphi\left(  t,s\right)  +C^{\top}M\left(
t,s\right)  \left(  \sigma-D\psi\left(  s\right)  \right) \\
\text{ \ \ \ \ \ \ }+%
{\displaystyle\int_{Z}}
E\left(  z\right)  ^{\top}M\left(  t,s\right)  \left(  c\left(  z\right)
-F\left(  z\right)  \psi\left(  s\right)  \right)  \theta\left(  dz\right)
,\\
M\left(  t,T\right)  =G\left(  t\right)  ,\text{ }\bar{M}\left(  t,T\right)
=\bar{G}\left(  t\right)  ,\text{\ }\Upsilon\left(  t,T\right)  =\mu
_{2}\left(  t\right)  \text{, }\varphi\left(  t,T\right)  =\mu_{2}\left(
t\right)  ,t\in\left[  0,T\right]  ,
\end{array}
\right.  \tag{4.11}%
\end{equation}
where
\[
\det\left[  R\left(  t,t\right)  +D\left(  t\right)  ^{\top}M\left(
t,t\right)  D\left(  t\right)  +\int_{Z}F\left(  t,z\right)  ^{\top}M\left(
t,t\right)  F\left(  t,z\right)  \theta\left(  dz\right)  \right]
\neq0,\text{ }t\in\left[  0,T\right]  ,
\]
and $\Psi\left(  s\right)  ,\psi\left(  s\right)  $ are given by $\left(
4.9\right)  .$

\subsection{Verification Theorem}

The following theorem provides a verification argument.

\begin{theorem}
let \textbf{(H1)-(H2) }hold$.$ Let $M\left(  .,.\right)  ,\bar{M}\left(
.,.\right)  ,\Upsilon\left(  .,.\right)  $ and $\varphi\left(  .,.\right)  $
be the solution of the system $\left(  4.11\right)  .$ Then $\hat{u}\left(
.\right)  $ given by $\left(  4.8\right)  $ is an equilibrium control.
\end{theorem}

\bop First, we can check that $\Psi\left(  .\right)  $ and $\psi\left(
.\right)  $\ in $\left(  3.15\right)  $ are both uniformly bounded. Then the
following linear SDE%
\[
\left\{
\begin{array}
[c]{l}%
d\hat{X}\left(  s\right)  =\left\{  \left(  A-B\Psi\left(  s\right)  \right)
\hat{X}\left(  s\right)  +b-B\psi\left(  s\right)  \right\}  ds\\
\text{ \ \ \ \ \ \ \ \ \ \ \ \ }+\left\{  \left(  C-D\Psi\left(  s\right)
\right)  \hat{X}\left(  s\right)  +\sigma-D\psi\left(  s\right)  \right\}
dW\left(  s\right) \\
\text{ \ \ \ \ \ \ \ \ \ \ \ \ }+%
{\displaystyle\int_{Z}}
\left\{  \left(  E\left(  z\right)  -F\left(  z\right)  \Psi\left(  s\right)
\right)  \hat{X}\left(  s-\right)  +c\left(  z\right)  -F\left(  z\right)
\psi\left(  s\right)  \right\}  \widetilde{N}\left(  ds,dz\right)  ,\text{ for
}s\in\left[  0,T\right]  ,\\
\hat{X}\left(  0\right)  =x_{0},
\end{array}
\right.
\]
is uniquely solvable, and the following estimate holds%
\[
\mathbb{E}\left[  \sup_{s\in\left[  0,T\right]  }\left\vert \hat{X}\left(
s\right)  \right\vert ^{2}\right]  \leq K\left(  1+x_{0}^{2}\right)  .
\]

So the control $\hat{u}\left(  .\right)  $ defined by $\left(  4.8\right)  $
is admissible. Moreover, by definition of $\left(  p\left(  s;t\right)
,q\left(  s;t\right)  ,r\left(  s,z;t\right)  \right)  $ via $\left(
4.3\right)  $, $\left(  4.6\right)  $ and $\left(  4.7\right)  ,$
respectively, and by applying the Ito's formula, we can easily check that, for
each $t\in\left[  0,T\right]  $ the processes $\left(  p\left(  .;t\right)
,q\left(  .;t\right)  ,r\left(  .,.;t\right)  \right)  $ satisfy $\left(
3.4\right)  .$

Finally, in view of Corollary $3.4,$ it's remains to check that the condition
$\left(  3.11\right)  $ holds. To this end, we substitute $\left(  p\left(
t;t\right)  ,q\left(  t;t\right)  ,r\left(  t,z;t\right)  ,\hat{u}\left(
t\right)  \right)  $ taken from $\left(  4.3\right)  ,$ $\left(  4.6\right)
,$ $\left(  4.7\right)  $ and $\left(  4.8\right)  ,$ respecively, in the
expression
\[
K\left(  t\right)  =R\left(  t,t\right)  \hat{u}\left(  t\right)  -B\left(
t\right)  ^{\top}p\left(  t;t\right)  -D\left(  t\right)  ^{\top}q\left(
t;t\right)  -%
{\displaystyle\int_{Z}}
F\left(  t,z\right)  ^{\top}r\left(  t,z;t\right)  \theta\left(  dz\right)  ,
\]
we get%
\begin{align*}
K\left(  t\right)   &  =R\left(  t,t\right)  \hat{u}\left(  t\right)
+B\left(  t\right)  ^{\top}\left\{  \left(  M\left(  t,t\right)  +\bar
{M}\left(  t,t\right)  +\Upsilon\left(  t,t\right)  \right)  \hat{X}\left(
t\right)  +\varphi\left(  t,t\right)  \right\}  +D\left(  t\right)  ^{\top
}M\left(  t,t\right)  \left(  C\left(  t\right)  \hat{X}\left(  t\right)
\right. \\
&  \text{ \ \ \ \ }+\left.  D\left(  t\right)  \hat{u}\left(  t\right)
+\sigma\left(  t\right)  \right)  +%
{\displaystyle\int_{Z}}
F\left(  t,z\right)  ^{\top}M\left(  t,t\right)  \left(  E\left(  t,z\right)
\hat{X}\left(  t\right)  +F\left(  t,z\right)  \hat{u}\left(  t\right)
+c\left(  t,z\right)  \right)  \theta\left(  dz\right) \\
&  =\Theta\left(  t\right)  ^{-1}\left(  \hat{u}\left(  t\right)  +\Psi\left(
t\right)  \hat{X}\left(  t\right)  +\psi\left(  t\right)  \right)  ,
\end{align*}
the representation $\left(  4.8\right)  $ shows that $K\left(  t\right)  =0$.
Hence, by Corollary $3.4,$ $\hat{u}\left(  .\right)  $ is an open-loop Nash
equilibrium control.\eop

Note that, the verification theorem (Theorem 4.1 ) assumes the existence of a
solution to the system $\left(  4.11\right)  $. However, since the ODEs which
should be solved by $M\left(  .,.\right)  $ and $\bar{M}\left(  .,.\right)  $
do not have a symmetry structure. The general solvability for this type of
ODEs when $(n>1)$ remains an outstanding open problem. We will see in the next
section two examples in the case when $n=1$, this case is important,
especially in financial applications as the mean--variance portfolio selection
model. Also, we remark that a special feature of the case when $n=1$ is that
the state $X\left(  .\right)  $ is one-dimensional, so are the unknown
$M\left(  .,.\right)  ,\bar{M}\left(  .,.\right)  ,\Upsilon\left(  .,.\right)
$ and $\varphi\left(  .,.\right)  $ of the system $\left(  4.11\right)  $.
This makes it easier to solve $\left(  4.11\right)  $.

\section{Some applications}

\subsection{Mean-variance portfolio selection problem}

In this subsection, we discuss the continuous-time Markowitz's mean-variance
portfolio selection problem. We apply Theorem 4.1 to obtain a state feedback
representation of an equilibrium control for the Problem$.$ In the absence of
Poisson random jumps this problem is discussed in \cite{18}$.$

The problem is formulated as follows: We consider a financial market, in which
two securities are traded continuously. One of them is a bond, with price
$S^{0}\left(  s\right)  $ at time $s\in\left[  0,T\right]  $ governed by
\begin{equation}
dS^{0}\left(  s\right)  =S^{0}\left(  s\right)  r\left(  s\right)  ds,\text{
}S^{0}\left(  0\right)  =s_{0}>0. \tag{5.1}%
\end{equation}

There is also a stock with unit price $S^{1}\left(  s\right)  $ at time
$s\in\left[  0,T\right]  $ governed by%
\begin{equation}
dS^{1}\left(  s\right)  =S^{1}\left(  s-\right)  \left(  \alpha\left(
s\right)  ds+\beta\left(  s\right)  dW\left(  s\right)  +\int_{%
\mathbb{R}
^{\ast}}\gamma\left(  s,z\right)  \widetilde{N}\left(  ds,dz\right)  \right)
,\text{ }S^{1}\left(  0\right)  =s^{1}>0. \tag{5.2}%
\end{equation}
where $r:\left[  0,T\right]  \rightarrow\left(  0,\infty\right)  ,$
$\alpha,\beta:\left[  0,T\right]  \rightarrow%
\mathbb{R}
$ and $\gamma:\left[  0,T\right]  \times%
\mathbb{R}
^{\ast}\rightarrow%
\mathbb{R}
$ are assumed to be deterministic and continuous$,$ we also assume a uniform
ellipticity condition as follow $\sigma\left(  t\right)  ^{2}+\int_{Z}%
\gamma\left(  t,z\right)  ^{2}\theta\left(  dz\right)  \geq\delta,$ $a.e$, for
some $\delta>0$. For an investor, a portfolio $\pi\left(  .\right)  $ is a
process represents the amount of money invested in the stock$.$ The wealth
process $X^{x_{0},\pi\left(  .\right)  }\left(  .\right)  $ corresponding to
initial capital $x_{0}>0,$ and portfolio $\pi\left(  .\right)  $, satisfies
then the equation%
\begin{equation}
\left\{
\begin{array}
[c]{l}%
dX\left(  s\right)  =\left(  r\left(  s\right)  X\left(  s\right)  +\pi\left(
s\right)  \left(  \alpha\left(  s\right)  -r\left(  s\right)  \right)
\right)  ds+\pi\left(  s\right)  \beta\left(  s\right)  dW\left(  s\right) \\
\text{ \ \ \ \ \ \ \ \ \ \ \ \ \ }+\pi\left(  s\right)
{\displaystyle\int_{\mathbb{R}^{\ast}}}
\gamma\left(  s,z\right)  \widetilde{N}\left(  ds,dz\right)  ,\text{ for }%
t\in\left[  0,T\right]  ,\\
X\left(  0\right)  =x_{0}.
\end{array}
\right.  \tag{5.3}%
\end{equation}

As time evolves, we need to consider the controlled stochastic differential
equation parametrized by $\left(  t,\xi\right)  \in\left[  0,T\right]
\times\mathbb{L}^{2}\left(  \Omega,\mathcal{F}_{t},\mathbb{P};%
\mathbb{R}
\right)  $ and satisfiied by $X\left(  .\right)  $%
\begin{equation}
\left\{
\begin{array}
[c]{l}%
dX\left(  s\right)  =\left(  r\left(  s\right)  X\left(  s\right)  +\pi\left(
s\right)  \left(  \alpha\left(  s\right)  -r\left(  s\right)  \right)
\right)  ds+\pi\left(  s\right)  \beta\left(  s\right)  dW\left(  s\right) \\
\text{ \ \ \ \ \ \ \ \ \ \ \ \ }+\pi\left(  s\right)
{\displaystyle\int_{\mathbb{R}^{\ast}}}
\gamma\left(  s,z\right)  \widetilde{N}\left(  ds,dz\right)  ,\text{ for }%
s\in\left[  t,T\right]  ,\\
X\left(  t\right)  =\xi.
\end{array}
\right.  \tag{5.4}%
\end{equation}

The objective is to maximize the conditional expectation of terminal wealth
$\mathbb{E}^{t}\left[  X\left(  T\right)  \right]  ,$ and at the same time to
minimize the conditional variance of the terminal wealth $\mathrm{Var}%
^{t}\left[  X\left(  T\right)  \right]  ,$ over controls $\pi\left(  .\right)
$ valued in $%
\mathbb{R}
$. Then, the mean-variance portfolio optimization problem is denoted as:
minimizing the cost $J\left(  t,\xi,.\right)  $, given by%
\begin{equation}
J\left(  t,\xi,\pi\left(  .\right)  \right)  =\frac{1}{2}\mathrm{Var}%
^{t}\left[  X\left(  T\right)  \right]  -\left(  \mu_{1}\left(  t\right)
\xi+\mu_{2}\left(  t\right)  \right)  \mathbb{E}^{t}\left[  X\left(  T\right)
\right]  , \tag{5.5}%
\end{equation}
subject to $\left(  5.4\right)  .$ Here $\mu_{1},\mu_{2}:\left[  0,T\right]
\rightarrow\left(  0,\infty\right)  ,$ are some deterministic non constant,
continuous and bounded functions. The above model cover the one in
\cite{18}$,$ since, in our case, the weight between the conditional variance
and the conditional expectation depends on the current wealth level, as well
as, the current time, while in \cite{18} the weight depends on the current
wealth level only. Hence, in the above model, there are three different
sources of time-inconsistency. Moreover, the above model is mathematically a
special case of the general LQ problem formulated earlier in this paper, with
$n=d=m=1$. Then we can apply Theorem 4.1 to obtain a Nash equilibrium control.
We recall that, the definition of equilibrium controls is in the sense of
open-loop, which is different from the feedback one in \cite{13}, \cite{19}
and \cite{9}.

The optimal control problem associated with $\left(  5.4\right)  $ and
$\left(  5.5\right)  $ is equivalent to minimize%
\[
J\left(  t,\xi,u\left(  .\right)  \right)  =\frac{1}{2}\left(  \mathbb{E}%
^{t}\left[  X\left(  T\right)  ^{2}\right]  -\mathbb{E}^{t}\left[  X\left(
T\right)  \right]  ^{2}\right)  -\left(  \mu_{1}\left(  t\right)  \xi+\mu
_{2}\left(  t\right)  \right)  \mathbb{E}^{t}\left[  X\left(  T\right)
\right]
\]
subject to $\left(  5.4\right)  .$ Denote%
\[
\rho\left(  t\right)  =\beta\left(  t\right)  ^{2}+\int_{%
\mathbb{R}
^{\ast}}\gamma\left(  t,z\right)  ^{2}\theta\left(  dz\right)  .
\]

Thus, the system $\left(  4.11\right)  $ reduces to
\begin{equation}
\left\{
\begin{array}
[c]{l}%
\dfrac{\partial M}{\partial s}\left(  t,s\right)  +\left\{  2r\left(
s\right)  -\dfrac{\left(  \alpha\left(  s\right)  -r\left(  s\right)  \right)
^{2}}{M\left(  s,s\right)  \rho\left(  s\right)  }\left(  M\left(  s,s\right)
+\bar{M}\left(  s,s\right)  +\Upsilon\left(  s,s\right)  \right)  \right\}
M\left(  t,s\right)  =0,\text{ }\forall s\in\left[  0,T\right]  ,\\
\dfrac{\partial\bar{M}}{\partial s}\left(  t,s\right)  +\left\{  2r\left(
s\right)  -\dfrac{\left(  \alpha\left(  s\right)  -r\left(  s\right)  \right)
^{2}}{M\left(  s,s\right)  \rho\left(  s\right)  }\left(  M\left(  s,s\right)
+\bar{M}\left(  s,s\right)  +\Upsilon\left(  s,s\right)  \right)  \right\}
\bar{M}\left(  t,s\right)  =0,\text{ }\forall s\in\left[  0,T\right]  ,\\
\dfrac{\partial\Upsilon}{\partial s}\left(  t,s\right)  +r\left(  s\right)
\Upsilon\left(  t,s\right)  =0,\text{ }\forall\left(  t,s\right)  \in D\left[
0,T\right]  ,\\
\dfrac{\partial\varphi}{\partial s}\left(  t,s\right)  +r\left(  s\right)
\varphi\left(  t,s\right)  =0,\text{ }\forall\left(  t,s\right)  \in D\left[
0,T\right]  ,\\
M\left(  t,T\right)  =1,\text{ }\bar{M}\left(  t,T\right)  =-1,\text{
}\Upsilon\left(  t,T\right)  =\mu_{1}\left(  t\right)  ,\text{ }\varphi\left(
t,T\right)  =\mu_{2}\left(  t\right)  ,\text{ }\forall t\in\left[  0,T\right]
.
\end{array}
\right.  \tag{5.6}%
\end{equation}

Clearly, if $M\left(  .,.\right)  $ and $\bar{M}\left(  .,.\right)  $ are
solutions to the first and the second equations, respectively, in $\left(
5.6\right)  $, then $\tilde{M}\left(  .,.\right)  =\left(  \bar{M}+M\right)
\left(  .,.\right)  $ solves the following ODE
\begin{equation}
\left\{
\begin{array}
[c]{l}%
\dfrac{\partial\tilde{M}}{\partial s}\left(  t,s\right)  +\left\{  2r\left(
s\right)  -\dfrac{\left(  \alpha\left(  s\right)  -r\left(  s\right)  \right)
^{2}}{M\left(  s,s\right)  \rho\left(  s\right)  }\left(  \tilde{M}\left(
s,s\right)  +\Upsilon\left(  s,s\right)  \right)  \right\}  \tilde{M}\left(
t,s\right)  =0,\text{ }\forall\left(  t,s\right)  \in\left[  0,T\right]  ,\\
\tilde{M}\left(  t,T\right)  =0,\text{ }t\in\left[  0,T\right]  ,
\end{array}
\right.  \tag{5.7}%
\end{equation}
which is equivalent to%
\begin{align*}
\tilde{M}\left(  t,s\right)   &  =\tilde{M}\left(  t,T\right)  e^{%
{\textstyle\int_{s}^{T}}
\left(  2r\left(  \tau\right)  -\tfrac{\left(  \alpha\left(  \tau\right)
-r\left(  \tau\right)  \right)  ^{2}}{M\left(  \tau,\tau\right)  \rho\left(
\tau\right)  }\left(  \tilde{M}\left(  \tau,\tau\right)  +\Upsilon\left(
\tau,\tau\right)  \right)  \right)  d\tau},\\
&  =0,\text{ }\left(  t,s\right)  \in D\left[  0,T\right]  ,
\end{align*}
hence%
\[
\bar{M}\left(  t,s\right)  +M\left(  t,s\right)  =\tilde{M}\left(  t,s\right)
=0,\text{ }\forall\left(  t,s\right)  \in D\left[  0,T\right]  .
\]

Moreover, we remark that all data of the ODEs which should be solved by
$M\left(  .,.\right)  $ and $\bar{M}\left(  .,.\right)  $ are not influenced
by $t,$ thus $\left(  5.6\right)  $ reduces to
\begin{equation}
\left\{
\begin{array}
[c]{l}%
\dfrac{dM}{ds}\left(  s\right)  +2r\left(  s\right)  M\left(  s\right)
-M\left(  s\right)  \dfrac{\left(  \alpha\left(  s\right)  -r\left(  s\right)
\right)  ^{2}}{M\left(  s\right)  \rho\left(  s\right)  }\Upsilon\left(
s,s\right)  =0,\text{ }\forall s\in\left[  0,T\right]  ,\\
\bar{M}\left(  s\right)  =-M\left(  s\right)  ,\text{ }\forall s\in\left[
0,T\right]  ,\\
\dfrac{\partial\Upsilon}{\partial s}\left(  t,s\right)  +r\left(  s\right)
\Upsilon\left(  t,s\right)  =0,\text{ }\forall\left(  t,s\right)  \in D\left[
0,T\right]  ,\\
\dfrac{\partial\varphi}{\partial s}\left(  t,s\right)  +r\left(  s\right)
\varphi\left(  t,s\right)  =0,\text{ }\forall\left(  t,s\right)  \in D\left[
0,T\right]  ,\\
M\left(  T\right)  =1,\text{ }\Upsilon\left(  t,T\right)  =\mu_{1}\left(
t\right)  ,\text{ }\varphi\left(  t,T\right)  =\mu_{2}\left(  t\right)
,\text{ }\forall t\in\left[  0,T\right]  .
\end{array}
\right.  \tag{5.8}%
\end{equation}
which is explicitly solved by%
\begin{equation}
\left\{
\begin{array}
[c]{l}%
M\left(  s\right)  =e^{2\int_{s}^{T}r\left(  \tau\right)  d\tau}\left\{  1+%
{\displaystyle\int_{s}^{T}}
e^{-\int_{\tau}^{T}r\left(  l\right)  dl}\mu_{1}\left(  \tau\right)
\dfrac{\left(  \alpha\left(  \tau\right)  -r\left(  \tau\right)  \right)
^{2}}{\rho\left(  \tau\right)  }d\tau\right\}  ,\text{ }\forall s\in\left[
0,T\right]  ,\\
\bar{M}\left(  s\right)  =-e^{2\int_{s}^{T}r\left(  \tau\right)  d\tau
}\left\{  1+%
{\displaystyle\int_{s}^{T}}
e^{-\int_{\tau}^{T}r\left(  l\right)  dl}\mu_{1}\left(  \tau\right)
\dfrac{\left(  \alpha\left(  \tau\right)  -r\left(  \tau\right)  \right)
^{2}}{\rho\left(  \tau\right)  }d\tau\right\}  ,\text{ }\forall s\in\left[
0,T\right]  ,\\
\Upsilon\left(  t,s\right)  =\mu_{1}\left(  t\right)  e^{%
{\textstyle\int_{s}^{T}}
r\left(  \tau\right)  d\tau},\text{ }\forall\left(  t,s\right)  \in D\left[
0,T\right]  ,\\
\varphi\left(  t,s\right)  =\mu_{2}\left(  t\right)  e^{%
{\textstyle\int_{s}^{T}}
r\left(  \tau\right)  d\tau},\text{ }\forall\left(  t,s\right)  \in D\left[
0,T\right]  .
\end{array}
\right.  \tag{5.9}%
\end{equation}

In view of Theorem 4.1, the representation of the Nash equilibrium\ control
$\left(  4.8\right)  $ then gives%
\begin{equation}
\hat{\pi}\left(  s\right)  =\Psi\left(  s\right)  \hat{X}\left(  s\right)
+\psi\left(  s\right)  ,\forall s\in\left[  0,T\right]  ,\nonumber
\end{equation}
where, $\forall s\in\left[  0,T\right]  $%
\begin{equation}
\Psi\left(  s\right)  =\dfrac{\left(  \alpha\left(  s\right)  -r\left(
s\right)  \right)  }{M\left(  s\right)  \rho\left(  s\right)  }\Upsilon\left(
s,s\right)  \text{ }\text{and }\psi\left(  s\right)  =\dfrac{\left(
\alpha\left(  s\right)  -r\left(  s\right)  \right)  }{M\left(  s\right)
\rho\left(  s\right)  }\varphi\left(  s,s\right)  .\nonumber
\end{equation}

The corresponding equilibrium dynamics solves the SDEJ%
\[
\left\{
\begin{array}
[c]{l}%
d\hat{X}\left(  s\right)  =\left\{  \left(  r\left(  s\right)  -\Psi\left(
s\right)  \left(  \alpha\left(  s\right)  -r\left(  s\right)  \right)
\right)  \hat{X}\left(  s\right)  -\psi\left(  s\right)  \left(  \alpha\left(
s\right)  -r\left(  s\right)  \right)  \right\}  ds\\
\text{ \ \ \ \ \ \ \ \ \ \ \ \ \ }-\left(  \Psi\left(  s\right)  \hat
{X}\left(  s\right)  +\psi\left(  s\right)  \right)  \left\{  \beta\left(
s\right)  dW\left(  s\right)  +%
{\displaystyle\int_{\mathbb{R}^{\ast}}}
\gamma\left(  s,z\right)  \widetilde{N}\left(  ds,dz\right)  \right\}  ,\text{
for }s\in\left[  0,T\right]  ,\\
\hat{X}\left(  0\right)  =x_{0}.
\end{array}
\right.
\]

\begin{remark}
In the absence of Poisson random jumps, we have the following items

\begin{enumerate}
\item In the case where $\mu_{1}\left(  t\right)  =\mu_{1}>0$ and $\mu
_{2}\left(  t\right)  =\mu_{2}>0,$ the solution obtained, for the
mean--variance problem, coincides with that obtained by \cite{18}.

\item In the case where $\mu_{1}\left(  t\right)  =\mu_{1}>0$ and $\mu
_{2}\left(  t\right)  =0,$ the solution obtained, coincides with that obtained
by \cite{13}.

\item In the case where $\mu_{1}\left(  t\right)  =0$ and $\mu_{2}\left(
t\right)  =\mu_{2}>0,$ the solution, obtained coincides with that obtained by
\cite{21}.
\end{enumerate}
\end{remark}

\subsection{General discounting LQ regulator}

In this subsection, we consider\ an example of a general discounting
time-inconsistent LQ model. The objective is to minimize the expected cost
functional, that is earned during a finite time horizon%
\begin{equation}
J\left(  t,\xi,u\left(  .\right)  \right)  =\frac{1}{2}\mathbb{E}^{t}\left[
\int_{t}^{T}\left\vert u\left(  s\right)  \right\vert ^{2}ds+h\left(
t\right)  \left\vert X\left(  T\right)  -\xi\right\vert ^{2}\right]
\tag{5.11}%
\end{equation}
where $h\left(  .\right)  :\left[  0,T\right]  \rightarrow\left(
0,\infty\right)  ,$ is a general deterministic non-exponential discount
function satisfying $h\left(  0\right)  =1,$ $h\left(  s\right)  \geq0$ and
$\int_{0}^{T}h\left(  t\right)  dt<\infty$. Subject to a controlled one
dimontional SDE parametrized by $\left(  t,\xi\right)  \in\left[  0,T\right]
\times\mathbb{L}^{2}\left(  \Omega,\mathcal{F}_{t},\mathbb{P};%
\mathbb{R}
\right)  $%
\begin{equation}
\left\{
\begin{array}
[c]{l}%
dX\left(  s\right)  =\left\{  aX\left(  s\right)  +bu\left(  s\right)
\right\}  ds+\sigma dW\left(  s\right)  +c%
{\displaystyle\int_{\mathbb{R}^{\ast}}}
\tilde{N}\left(  ds,dz\right)  ,\text{ }s\in\left[  0,T\right]  ,\\
X\left(  t\right)  =\xi,
\end{array}
\right.  \tag{5.12}%
\end{equation}
where $a,b$ and $c$ are real constant$.$ As mentioned in \cite{8}$,$ this is
atime-inconsistent version of the classical linear quadratic regulator, we
want control the system so that the final sate $X\left(  T\right)  $ is close
to $\xi$ while at the same time we keep the control energy (formalized by the
integral term) small. Note that, here the time-inconsistency is due to the
fact that the terminal cost depends explicitly on the current state $\xi$ as
well as the current time $t.$ Hence there are two different sources of
time-inconsistency. For this example, the system $\left(  4.11\right)  $
reduces to%
\begin{equation}
\left\{
\begin{array}
[c]{l}%
\dfrac{\partial M}{\partial s}\left(  t,s\right)  +2aM\left(  t,s\right)
-b^{2}M\left(  t,s\right)  \left\{  M\left(  s,s\right)  +\bar{M}\left(
s,s\right)  +\Upsilon\left(  s,s\right)  \right\}  =0,\text{ }\forall\left(
t,s\right)  \in D\left[  0,T\right]  ,\\
\dfrac{\partial\bar{M}}{\partial s}\left(  t,s\right)  +2a\bar{M}\left(
t,s\right)  -b^{2}\bar{M}\left(  t,s\right)  \left\{  M\left(  s,s\right)
+\bar{M}\left(  s,s\right)  +\Upsilon\left(  s,s\right)  \right\}  =0,\text{
}\forall\left(  t,s\right)  \in D\left[  0,T\right]  ,\\
\dfrac{\partial\Upsilon}{\partial s}\left(  t,s\right)  +a\Upsilon\left(
t,s\right)  =0,\text{ }\forall\left(  t,s\right)  \in D\left[  0,T\right]  ,\\
\dfrac{\partial\varphi}{\partial s}\left(  t,s\right)  +a\varphi\left(
t,s\right)  -b^{2}\left\{  M\left(  t,s\right)  +\bar{M}\left(  t,s\right)
\right\}  \varphi\left(  s,s\right)  =0,\text{ }\forall\left(  t,s\right)  \in
D\left[  0,T\right]  ,\\
M\left(  t,T\right)  =h\left(  t\right)  ,\text{ }\bar{M}\left(  t,T\right)
=0,\text{\ }\Upsilon\left(  t,T\right)  =h\left(  t\right)  \text{, }%
\varphi\left(  t,T\right)  =0,\text{ }\forall t\in\left[  0,T\right]  ,
\end{array}
\right.  \tag{5.13}%
\end{equation}
obviously, $\Upsilon\left(  .,.\right)  $ is explicitely given by%
\begin{equation}
\Upsilon\left(  t,s\right)  =h\left(  t\right)  \exp\left\{  a\left(
T-s\right)  \right\}  ,\text{ }\forall\left(  t,s\right)  \in D\left[
0,T\right]  . \tag{5.14}%
\end{equation}

Moreover, we can check that $M\left(  .,.\right)  ,\bar{M}\left(  .,.\right)
$ and $\varphi\left(  .,.\right)  $ solve $\left(  5.13\right)  ,$ if and only
if, they solve the following system of integral equations
\begin{equation}
\left\{
\begin{array}
[c]{l}%
M\left(  t,s\right)  =M\left(  t,T\right)  e^{%
{\textstyle\int_{s}^{T}}
\left\{  2a-b^{2}\left(  M\left(  \tau,\tau\right)  +\bar{M}\left(  \tau
,\tau\right)  +\Upsilon\left(  \tau,\tau\right)  \right)  \right\}  d\tau
}\text{, }\forall\left(  t,s\right)  \in D\left[  0,T\right]  .\\
\bar{M}\left(  t,s\right)  =\bar{M}\left(  t,T\right)  e^{%
{\textstyle\int_{s}^{T}}
\left\{  2a-b^{2}\left(  M\left(  \tau,\tau\right)  +\bar{M}\left(  \tau
,\tau\right)  +\Upsilon\left(  \tau,\tau\right)  \right)  \right\}  d\tau
},\text{\ }\forall\left(  t,s\right)  \in D\left[  0,T\right]  .\\
\varphi\left(  t,s\right)  =\varphi\left(  t,T\right)  e^{a\left(  T-s\right)
}-b^{2}%
{\displaystyle\int_{s}^{T}}
e^{a\left(  \tau-s\right)  }\left(  M\left(  t,\tau\right)  +\bar{M}\left(
t,\tau\right)  \right)  \varphi\left(  \tau,\tau\right)  d\tau,\text{\ }%
\forall\left(  t,s\right)  \in D\left[  0,T\right]  ,
\end{array}
\right.  \tag{5.15}%
\end{equation}
on the other hand, we have $\bar{M}\left(  t,T\right)  =\varphi\left(
t,T\right)  =0,$ then $\left(  5.15\right)  $ reduces to
\begin{equation}
\left\{
\begin{array}
[c]{l}%
M\left(  t,s\right)  =M\left(  t,T\right)  e^{%
{\textstyle\int_{s}^{T}}
\left\{  2a-b^{2}\left(  M\left(  r,r\right)  +\Upsilon\left(  r,r\right)
\right)  \right\}  dr},\text{ }\forall\left(  t,s\right)  \in D\left[
0,T\right]  .\\
\bar{M}\left(  t,s\right)  =0,\text{ }\forall\left(  t,s\right)  \in D\left[
0,T\right]  ,\\
\varphi\left(  t,s\right)  =-b^{2}%
{\displaystyle\int_{s}^{T}}
e^{a\left(  \tau-s\right)  }M\left(  t,\tau\right)  \varphi\left(  \tau
,\tau\right)  d\tau,\text{ }\forall\left(  t,s\right)  \in D\left[
0,T\right]  .
\end{array}
\right.  \tag{5.16}%
\end{equation}

It is clear that if $M\left(  .,.\right)  $ is the solution of the first
equation in $\left(  5,16\right)  ,$ then%
\[
\varphi\left(  s,s\right)  =-b^{2}%
{\displaystyle\int_{s}^{T}}
e^{a\left(  \tau-s\right)  }M\left(  s,\tau\right)  \varphi\left(  \tau
,\tau\right)  d\tau,\text{ }\forall s\in\left[  0,T\right]  ,
\]
thus, there exists some constant $L>0$ such that $\left\vert \varphi\left(
s,s\right)  \right\vert \leq L%
{\displaystyle\int_{s}^{T}}
\left\vert \varphi\left(  \tau,\tau\right)  \right\vert d\tau$, then by
Gronwall Lemma, we conclude that $\varphi\left(  s,s\right)  =0,\forall
s\in\left[  0,T\right]  .$ Therefore%
\[
\varphi\left(  t,s\right)  =0,\text{ }\forall\left(  t,s\right)  \in D\left[
0,T\right]  .
\]
is the unique solution to the last equation in the system $\left(
5.16\right)  .$

Now, it's remains to solve the first equation in the system $\left(
5.16\right)  .$ We have the following Theorem.

\begin{theorem}
The first equation in $\left(  5.16\right)  $ has a unique solution in
$C\left(  D\left[  0,T\right]  ,%
\mathbb{R}
^{+}\right)  .$
\end{theorem}

\bop See the proof in Appendix \textbf{A.2.}\eop

In view of Theorem 4.1, the representation $\left(  4.8\right)  $ of the Nash
equilibrium control, then gives%
\begin{equation}
\hat{u}\left(  s\right)  =b\left\{  \Upsilon\left(  s,s\right)  -M\left(
s,s\right)  \right\}  \hat{X}\left(  s\right)  ,\forall s\in\left[
0,T\right]  , \tag{5.19}%
\end{equation}
and the corresponding equilibrium dynamics solves the SDEJ%
\begin{equation}
\left\{
\begin{array}
[c]{l}%
d\hat{X}\left(  s\right)  =\left\{  a+b^{2}\left(  \Upsilon\left(  s,s\right)
-M\left(  s,s\right)  \right)  \right\}  \hat{X}\left(  s\right)  ds+\sigma
dW\left(  s\right)  +c%
{\displaystyle\int_{\mathbb{R}^{\ast}}}
z\tilde{N}\left(  ds,dz\right)  ,\text{ }s\in\left[  0,T\right]  ,\\
X\left(  0\right)  =x_{0}.
\end{array}
\right.  \tag{5.20}%
\end{equation}

To conclude this section let us present the following remark.

\begin{remark}
The Problem (E) given by the subsection 2.3, is in fact shown to be a
particular case of the general discounting LQ regulator model, formulated
earlier in this paragraph, in the case when $a=c=0,$ and the final data
$\xi=x$, this leads to the following representation of the Nash equilibrium
control of this problem
\[
\hat{u}\left(  s\right)  =b\left(  h\left(  s\right)  -M\left(  s,s\right)
\right)  \hat{X}\left(  s\right)  ,\forall s\in\left[  0,T\right]  ,
\]
where $M\left(  t,s\right)  $ solves%
\[
M\left(  t,s\right)  =h\left(  t\right)  e^{%
{\textstyle\int_{s}^{T}}
-b^{2}\left(  M\left(  \tau,\tau\right)  +h\left(  \tau\right)  \right)
d\tau},\text{ for }\left(  t,s\right)  \in D\left[  0,T\right]  ,
\]
and the corresponding equilibrium dynamics solves the SDE%
\[
\left\{
\begin{array}
[c]{l}%
d\hat{X}\left(  s\right)  =b^{2}\left\{  h(s)-M\left(  s,s\right)  \right\}
\hat{X}\left(  s\right)  ds+\sigma dW\left(  s\right)  ,\text{ }s\in\left[
0,T\right]  ,\\
X\left(  0\right)  =x_{0}.
\end{array}
\right.
\]
This in fact, the correct solution of Problem (E). \bigskip
\end{remark}

{\Large Conclusion and future work.}\textit{ In this paper, we have studied a
class of dynamic decision problems of a general time-inconsistent type. We
have used the game theoretic approach to handle the time inconsistency. During
this study open-loop Nash equilibrium controls are constructed as an
alternative of optimal controls. This has been accomplished through stochastic
maximum principle that includes a flow of forward-backward stochastic
differential equations under maximum condition. The inclusion of concrete
examples confirms the validity of our proposed study. The work can be extended
in several ways. For example, this approach can be extended to a mean field
game to construct decentralized strategies and obtain an estimate of their
performance. The reserch on this topic is in progress and will appear in our
forthcoming paper.}

\section{Appendix: Additional proofs}

\textbf{A.1.} \textbf{Proof of Lamma 3.3.} Let $t\in\left[  0,T\right]  ,$
$v\in\mathbb{L}^{2}\left(  \Omega,\mathcal{F}_{t},\mathbb{P};%
\mathbb{R}
^{m}\right)  $ and $\varepsilon\in\left[  0,T-t\right]  .$ Since
$\mathbb{E}^{t}\left[  y^{\varepsilon,v}\left(  .\right)  \right]  $ and
$\mathbb{E}^{t}\left[  z^{\varepsilon,v}\left(  .\right)  \right]  $ solve the
following ODEs, respectively%
\[
\left\{
\begin{array}
[c]{l}%
d\mathbb{E}^{t}\left[  y^{\varepsilon,v}\left(  s\right)  \right]
=A\mathbb{E}^{t}\left[  y^{\varepsilon,v}\left(  s\right)  \right]  ds,\text{
}s\in\left[  t,T\right]  ,\\
\mathbb{E}^{t}\left[  y^{\varepsilon,v}\left(  t\right)  \right]  =0,\text{\ }%
\end{array}
\right.
\]
and%
\[
\left\{
\begin{array}
[c]{l}%
d\mathbb{E}^{t}\left[  z^{\varepsilon,v}\left(  s\right)  \right]  =\left\{
A\mathbb{E}^{t}\left[  z^{\varepsilon,v}\left(  s\right)  \right]
+B\mathbb{E}^{t}\left[  v\right]  1_{\left[  t,t+\varepsilon\right)  }\left(
s\right)  \right\}  ds,\text{ }s\in\left[  t,T\right]  ,\\
\mathbb{E}^{t}\left[  z^{\varepsilon,v}\left(  t\right)  \right]  =0.
\end{array}
\right.
\]

Thus, it is clear that $\mathbb{E}^{t}\left[  y^{\varepsilon,v}\left(
s\right)  \right]  =0,$ $a.e.$ $s\in\left[  t,T\right]  .$ According to
Gronwall's inequality there exists a positive constant $K$ such that
$\sup_{s\in\left[  t,T\right]  }\left\vert \mathbb{E}^{t}\left[
z^{\varepsilon,v}\left(  s\right)  \right]  \right\vert ^{2}\leq
K\varepsilon^{2}.$ Moreover, by Lemma 2.1. in \cite{15}, we obtain (3.18).

By these estimates, we can calculate the difference%
\begin{equation}%
\begin{array}
[c]{l}%
J\left(  t,\hat{X}\left(  t\right)  ,u^{\varepsilon}\left(  .\right)  \right)
-J\left(  t,\hat{X}\left(  t\right)  ,\hat{u}\left(  .\right)  \right) \\
=\mathbb{E}^{t}\left[
{\displaystyle\int_{t}^{T}}
\left\{  \left\langle Q\left(  t,s\right)  \hat{X}\left(  s\right)  +\bar
{Q}\left(  t,s\right)  \mathbb{E}^{t}\left[  \hat{X}\left(  s\right)  \right]
,y^{\varepsilon,v}\left(  s\right)  +z^{\varepsilon,v}\left(  s\right)
\right\rangle \right.  \right. \\
\text{ \ \ \ \ \ \ \ \ \ \ \ \ \ \ }+\dfrac{1}{2}\left\langle Q\left(
t,s\right)  \left(  y^{\varepsilon,v}\left(  s\right)  +z^{\varepsilon
,v}\left(  s\right)  \right)  ,y^{\varepsilon,v}\left(  s\right)
+z^{\varepsilon,v}\left(  s\right)  \right\rangle \\
\text{ \ \ \ \ \ \ \ \ \ \ \ \ \ \ }+\dfrac{1}{2}\left\langle \bar{Q}\left(
t,s\right)  \mathbb{E}^{t}\left[  y^{\varepsilon,v}\left(  s\right)
+z^{\varepsilon,v}\left(  s\right)  \right]  ,\mathbb{E}^{t}\left[
y^{\varepsilon,v}\left(  s\right)  +z^{\varepsilon,v}\left(  s\right)
\right]  \right\rangle \\
\text{ \ \ \ \ \ \ \ \ \ \ \ \ \ \ }+\left\langle R\left(  t,s\right)  \hat
{u}\left(  s\right)  ,v\right\rangle 1_{\left[  t,t+\varepsilon\right)
}\left(  s\right)  +\left.  \dfrac{1}{2}\left\langle R\left(  t,s\right)
v,v\right\rangle 1_{\left[  t,t+\varepsilon\right)  }\left(  s\right)
\right\}  ds\\
\text{ \ \ \ \ \ \ \ \ }+\dfrac{1}{2}\left\langle G\left(  t\right)  \left(
y^{\varepsilon,v}\left(  T\right)  +z^{\varepsilon,v}\left(  T\right)
\right)  ,y^{\varepsilon,v}\left(  T\right)  +z^{\varepsilon,v}\left(
T\right)  \right\rangle \\
\text{ \ \ \ \ \ \ \ \ }+\left\langle G\left(  t\right)  \hat{X}\left(
T\right)  +\bar{G}\left(  t\right)  \mathbb{E}^{t}\left[  \hat{X}\left(
T\right)  \right]  +\mu_{1}\left(  t\right)  \hat{X}\left(  t\right)  +\mu
_{2}\left(  t\right)  ,y^{\varepsilon,v}\left(  T\right)  +z^{\varepsilon
,v}\left(  T\right)  \right\rangle \\
\text{ \ \ \ \ \ \ \ \ }+\left.  \dfrac{1}{2}\left\langle \bar{G}\left(
t\right)  \mathbb{E}^{t}\left[  y^{\varepsilon,v}\left(  T\right)
+z^{\varepsilon,v}\left(  T\right)  \right]  ,\mathbb{E}^{t}\left[
y^{\varepsilon,v}\left(  T\right)  +z^{\varepsilon,v}\left(  T\right)
\right]  \right\rangle \right]  .
\end{array}
\tag{6.1}%
\end{equation}

In an other hand, from \textbf{(H1)} and $\left(  3.17\right)  -\left(
3.18\right)  $ the following estimate holds%
\begin{align*}
\mathbb{E}^{t}\left[
{\displaystyle\int_{t}^{T}}
\dfrac{1}{2}\left\langle \bar{Q}\left(  t,s\right)  \mathbb{E}^{t}\left[
y^{\varepsilon,v}\left(  s\right)  +z^{\varepsilon,v}\left(  s\right)
\right]  ,\mathbb{E}^{t}\left[  y^{\varepsilon,v}\left(  s\right)
+z^{\varepsilon,v}\left(  s\right)  \right]  \right\rangle ds\right.   & \\
+\left.  \dfrac{1}{2}\left\langle \bar{G}\left(  t\right)  \mathbb{E}%
^{t}\left[  y^{\varepsilon,v}\left(  T\right)  +z^{\varepsilon,v}\left(
T\right)  \right]  ,\mathbb{E}^{t}\left[  y^{\varepsilon,v}\left(  T\right)
+z^{\varepsilon,v}\left(  T\right)  \right]  \right\rangle \right]   &
=o\left(  \varepsilon\right)  .
\end{align*}

Then, from the terminal conditions in the adjoint equations, it follows that%
\begin{equation}%
\begin{array}
[c]{l}%
J\left(  t,\hat{X}\left(  t\right)  ,u^{\varepsilon}\left(  .\right)  \right)
-J\left(  t,\hat{X}\left(  t\right)  ,\hat{u}\left(  .\right)  \right) \\
=\mathbb{E}^{t}\left[
{\displaystyle\int_{t}^{T}}
\left\{  \left\langle Q\left(  t,s\right)  \hat{X}\left(  s\right)  +\bar
{Q}\left(  t,s\right)  \mathbb{E}^{t}\left[  \hat{X}\left(  s\right)  \right]
,y^{\varepsilon,v}\left(  s\right)  +z^{\varepsilon,v}\left(  s\right)
\right\rangle \right.  \right. \\
\text{ \ \ \ \ \ \ \ \ \ \ \ \ \ \ \ }+\dfrac{1}{2}\left.  \left\langle
Q\left(  t,s\right)  \left(  y^{\varepsilon,v}\left(  s\right)
+z^{\varepsilon,v}\left(  s\right)  \right)  ,y^{\varepsilon,v}\left(
s\right)  +z^{\varepsilon,v}\left(  s\right)  \right\rangle \right. \\
\text{ \ \ \ \ \ \ \ \ \ \ \ \ \ \ \ }+\left\langle R\left(  t,s\right)
\hat{u}\left(  s\right)  ,v\right\rangle 1_{\left[  t,t+\varepsilon\right)
}\left(  s\right)  +\left.  \dfrac{1}{2}\left\langle R\left(  t,s\right)
v,v\right\rangle 1_{\left[  t,t+\varepsilon\right)  }\left(  s\right)
\right\}  ds\\
\text{ \ \ \ \ \ \ \ \ }-\left\langle p\left(  T;t\right)  ,y^{\varepsilon
,v}\left(  T\right)  +z^{\varepsilon,v}\left(  T\right)  \right\rangle
-\left.  \dfrac{1}{2}\left\langle P\left(  T;t\right)  \left(  y^{\varepsilon
,v}\left(  T\right)  +z^{\varepsilon,v}\left(  T\right)  \right)
,y^{\varepsilon,v}\left(  T\right)  +z^{\varepsilon,v}\left(  T\right)
\right\rangle \right] \\
\text{ \ \ \ \ \ \ \ \ }+o\left(  \varepsilon\right)  .
\end{array}
\tag{6.2}%
\end{equation}

Now, by applying Ito's formula to $s\mapsto\left\langle p\left(  s;t\right)
,y^{\varepsilon,v}\left(  s\right)  +z^{\varepsilon,v}\left(  s\right)
\right\rangle $ on $\left[  t,T\right]  $, we get%
\begin{equation}%
\begin{array}
[c]{l}%
\left\langle p\left(  T;t\right)  ,y^{\varepsilon,v}\left(  T\right)
+z^{\varepsilon,v}\left(  T\right)  \right\rangle \\
=%
{\displaystyle\int\limits_{t}^{T}}
\left\{  \left(  Bv\right)  ^{\top}p\left(  s;t\right)  1_{\left[
t,t+\varepsilon\right)  }\left(  s\right)  \right.  +\left(  y^{\varepsilon
,v}\left(  s\right)  +z^{\varepsilon,v}\left(  s\right)  \right)  ^{\top
}\left(  Q\left(  t,s\right)  \hat{X}\left(  s\right)  +\bar{Q}\left(
t,s\right)  \mathbb{E}^{t}\left[  \hat{X}\left(  s\right)  \right]  \right) \\
\text{ \ \ \ \ \ }+%
{\displaystyle\sum\limits_{j=1}^{d}}
\left(  D_{j}v\right)  ^{\top}q_{j}\left(  s;t\right)  1_{\left[
t,t+\varepsilon\right)  }\left(  s\right)  +\left.
{\displaystyle\int_{Z}}
\left(  F\left(  z\right)  v\right)  ^{\top}r\left(  s,z;t\right)  1_{\left[
t,t+\varepsilon\right)  }\left(  s\right)  \theta\left(  dz\right)  \right\}
ds\\
+%
{\displaystyle\sum\limits_{j=1}^{d}}
{\displaystyle\int\limits_{t}^{T}}
\left\{  \left(  C_{j}\left(  y^{\varepsilon,v}\left(  s\right)
+z^{\varepsilon,v}\left(  s\right)  \right)  +D_{j}v1_{\left[  t,t+\varepsilon
\right)  }\left(  s\right)  \right)  ^{\top}p\left(  s;t\right)  +\left(
y^{\varepsilon,v}\left(  s\right)  +z^{\varepsilon,v}\left(  s\right)
\right)  ^{\top}q_{j}\left(  s;t\right)  \right\}  dW^{j}\left(  s\right) \\
+%
{\displaystyle\int\limits_{t}^{T}}
{\displaystyle\int\limits_{Z}}
\left\{  \left(  E\left(  z\right)  \left(  y^{\varepsilon,v}\left(  s\right)
+z^{\varepsilon,v}\left(  s\right)  \right)  +F\left(  z\right)  v1_{\left[
t,t+\varepsilon\right)  }\left(  s\right)  \right)  ^{\top}p\left(
s;t\right)  \right. \\
\text{ \ \ \ \ \ \ \ \ \ \ \ \ \ \ }+\left.  \left(  y^{\varepsilon,v}\left(
s\right)  +z^{\varepsilon,v}\left(  s\right)  \right)  ^{\top}r\left(
s,z;t\right)  \right\}  \tilde{N}\left(  ds,dz\right)  .
\end{array}
\tag{6.3}%
\end{equation}

By applying Ito's formula to $s\mapsto\left\langle P\left(  s;t\right)
\left(  y^{\varepsilon,v}\left(  s\right)  +z^{\varepsilon,v}\left(  s\right)
\right)  ,y^{\varepsilon,v}\left(  s\right)  +z^{\varepsilon,v}\left(
s\right)  \right\rangle $ on $\left[  t,T\right]  ,$ we get%
\begin{equation}%
\begin{array}
[c]{l}%
\left\langle P\left(  T;t\right)  \left(  y^{\varepsilon,v}\left(  T\right)
+z^{\varepsilon,v}\left(  T\right)  \right)  ,y^{\varepsilon,v}\left(
T\right)  +z^{\varepsilon,v}\left(  T\right)  \right\rangle \\
=%
{\displaystyle\int\limits_{t}^{T}}
\left\{  2\left(  y^{\varepsilon,v}\left(  s\right)  +z^{\varepsilon,v}\left(
s\right)  \right)  ^{\top}P\left(  s;t\right)  Bv1_{\left[  t,t+\varepsilon
\right)  }\left(  s\right)  +\left(  y^{\varepsilon,v}\left(  s\right)
+z^{\varepsilon,v}\left(  s\right)  \right)  ^{\top}Q\left(  t,s\right)
\left(  y^{\varepsilon,v}\left(  s\right)  +z^{\varepsilon,v}\left(  s\right)
\right)  \right. \\
\text{ \ \ \ \ \ \ \ \ }+%
{\displaystyle\sum\limits_{j=1}^{d}}
\left(  2\left(  y^{\varepsilon,v}\left(  s\right)  +z^{\varepsilon,v}\left(
s\right)  \right)  ^{\top}C_{j}^{\top}P\left(  s;t\right)  D_{j}v1_{\left[
t,t+\varepsilon\right)  }\left(  s\right)  \right.  \medskip\\
\text{ \ \ \ \ \ \ \ \ \ \ \ \ \ \ \ \ \ \ \ }+\left.  v^{\top}D_{j}^{\top
}P\left(  s;t\right)  D_{j}v1_{\left[  t,t+\varepsilon\right)  }\left(
s\right)  \right)  \medskip\\
\text{ \ \ \ \ \ \ }+%
{\displaystyle\int\limits_{Z}}
\left\{  2\left(  y^{\varepsilon,v}\left(  s\right)  +z^{\varepsilon,v}\left(
s\right)  \right)  ^{\top}E\left(  z\right)  ^{\top}P\left(  s;t\right)
F\left(  z\right)  v1_{\left[  t,t+\varepsilon\right)  }\left(  s\right)
\right. \\
\text{ \ \ \ \ \ \ \ \ \ \ \ \ \ \ \ }+\left.  v^{\top}F\left(  z\right)
^{\top}P\left(  s;t\right)  F\left(  z\right)  v1_{\left[  t,t+\varepsilon
\right)  }\left(  s\right)  \theta\left(  dz\right)  \right\}  ds\\
\text{ \ \ }+2%
{\displaystyle\sum\limits_{j=1}^{d}}
{\displaystyle\int\limits_{t}^{T}}
\left\{  \left(  y^{\varepsilon,v}\left(  s\right)  +z^{\varepsilon,v}\left(
s\right)  \right)  ^{\top}P\left(  s;t\right)  \left(  C_{j}\left(
y^{\varepsilon,v}\left(  s\right)  +z^{\varepsilon,v}\left(  s\right)
\right)  +D_{j}v1_{\left[  t,t+\varepsilon\right)  }\left(  s\right)  \right)
\right\}  dW^{j}\left(  s\right) \\
+2%
{\displaystyle\int\limits_{t}^{T}}
{\displaystyle\int\limits_{Z}}
\left\{  \left(  y^{\varepsilon,v}\left(  s\right)  +z^{\varepsilon,v}\left(
s\right)  \right)  ^{\top}P\left(  s;t\right)  \left(  E\left(  z\right)
\left(  y^{\varepsilon,v}\left(  s\right)  +z^{\varepsilon,v}\left(  s\right)
\right)  +F\left(  z\right)  v1_{\left[  t,t+\varepsilon\right)  }\left(
s\right)  \right)  \right\}  \tilde{N}\left(  ds,dz\right)  ,
\end{array}
\tag{6.4}%
\end{equation}

Moreover, we conclude from \textbf{(H1)} together with $\left(  3.17\right)
-\left(  3.18\right)  $ that%
\begin{equation}%
\begin{array}
[c]{l}%
\mathbb{E}^{t}\left[
{\displaystyle\int_{t}^{T}}
\left(  y^{\varepsilon,v}\left(  s\right)  +z^{\varepsilon,v}\left(  s\right)
\right)  ^{\top}P\left(  s;t\right)  Bv1_{\left[  t,t+\varepsilon\right)
}\left(  s\right)  ds\right]  =o\left(  \varepsilon\right)  ,\\
\mathbb{E}^{t}\left[
{\displaystyle\int_{t}^{T}}
\left(  y^{\varepsilon,v}\left(  s\right)  +z^{\varepsilon,v}\left(  s\right)
\right)  ^{\top}C_{j}^{\top}P\left(  s;t\right)  D_{j}v1_{\left[
t,t+\varepsilon\right)  }\left(  s\right)  ds\right]  =o\left(  \varepsilon
\right)  ,\\
\mathbb{E}^{t}\left[
{\displaystyle\int_{t}^{T}}
{\displaystyle\int_{Z}}
\left(  y^{\varepsilon,v}\left(  s\right)  +z^{\varepsilon,v}\left(  s\right)
\right)  ^{\top}E\left(  z\right)  ^{\top}P\left(  s;t\right)  F\left(
z\right)  v1_{\left[  t,t+\varepsilon\right)  }\left(  s\right)  \theta\left(
dz\right)  ds\right]  =o\left(  \varepsilon\right)  .
\end{array}
\tag{6.5}%
\end{equation}

By taking the conditional expectation in $\left(  6.3\right)  $ and $\left(
6.4\right)  ,$ then by invoking $\left(  6.5\right)  $ it hold that%
\begin{align}
&  \mathbb{E}^{t}\left[  \left\langle p\left(  T;t\right)  ,y^{\varepsilon
,v}\left(  T\right)  +z^{\varepsilon,v}\left(  T\right)  \right\rangle \right]
\nonumber\\
&  =\mathbb{E}^{t}\left[  \int\limits_{t}^{T}\left\{  v^{T}B^{\top}p\left(
s;t\right)  1_{\left[  t,t+\varepsilon\right)  }\left(  s\right)  \right.
\right.  +\left(  y^{\varepsilon,v}\left(  s\right)  +z^{\varepsilon,v}\left(
s\right)  \right)  ^{\top}\left(  Q\left(  t,s\right)  \hat{X}\left(
s\right)  +\bar{Q}\left(  t,s\right)  \mathbb{E}^{t}\left[  \hat{X}\left(
s\right)  \right]  \right) \nonumber\\
&  \text{ \ \ \ \ \ \ \ \ \ \ \ \ \ }+\sum\limits_{j=1}^{d}v^{T}D_{j}^{\top
}q_{j}\left(  s;t\right)  1_{\left[  t,t+\varepsilon\right)  }\left(
s\right)  +\left.  \left.
{\displaystyle\int_{Z}}
v^{T}F\left(  z\right)  ^{\top}r\left(  s,z;t\right)  1_{\left[
t,t+\varepsilon\right)  }\left(  s\right)  \theta\left(  dz\right)  \right\}
ds\right]  , \tag{6.6}%
\end{align}
and%
\begin{align}
&  \frac{1}{2}\mathbb{E}^{t}\left[  \left\langle P\left(  T;t\right)  \left(
y^{\varepsilon,v}\left(  T\right)  +z^{\varepsilon,v}\left(  T\right)
\right)  ,y^{\varepsilon,v}\left(  T\right)  +z^{\varepsilon,v}\left(
T\right)  \right\rangle \right] \nonumber\\
&  =\frac{1}{2}\mathbb{E}^{t}\left[  \int\limits_{t}^{T}\left\{  \left(
y^{\varepsilon,v}\left(  s\right)  +z^{\varepsilon,v}\left(  s\right)
\right)  ^{\top}Q\left(  t,s\right)  \left(  y^{\varepsilon,v}\left(
s\right)  +z^{\varepsilon,v}\left(  s\right)  \right)  \right.  \right.
+\sum\limits_{j=1}^{d}v^{\top}D_{j}^{\top}P\left(  s;t\right)  D_{j}%
v1_{\left[  t,t+\varepsilon\right)  }\left(  s\right) \nonumber\\
&  \text{ \ \ \ \ \ \ \ \ \ \ \ \ \ \ \ \ \ \ }+\text{\ }\left.  \left.
{\displaystyle\int_{Z}}
v^{\top}F\left(  z\right)  ^{\top}P\left(  s;t\right)  F\left(  z\right)
v1_{\left[  t,t+\varepsilon\right)  }\left(  s\right)  \theta\left(
dz\right)  \right\}  ds\right]  +o\left(  \varepsilon\right)  . \tag{6.7}%
\end{align}

By taking $\left(  6.6\right)  $ and $\left(  6.7\right)  $ in $\left(
6.2\right)  ,$ it follows that%
\begin{align*}
&  J\left(  t,\hat{X}\left(  t\right)  ,u^{\varepsilon}\left(  .\right)
\right)  -J\left(  t,\hat{X}\left(  t\right)  ,\hat{u}\left(  .\right)
\right) \\
&  =-\mathbb{E}^{t}\left[
{\displaystyle\int_{t}^{T}}
\left\{  v^{\top}B^{\top}p\left(  s;t\right)  +\sum\limits_{j=1}^{d}v^{\top
}D_{j}^{\top}q_{j}\left(  s,t\right)  \right.  +\frac{1}{2}\sum\limits_{j=1}%
^{d}v^{\top}D_{j}^{\top}P\left(  s;t\right)  D_{j}v\right. \\
&  \text{ \ \ \ \ \ \ \ \ \ \ \ \ \ \ \ \ \ \ \ \ \ }\left.  -v^{\top}R\left(
t,s\right)  \hat{u}\left(  s\right)  -\dfrac{1}{2}v^{\top}R\left(  t,s\right)
v\right. \\
&  \text{ \ \ \ \ \ \ \ \ \ \ \ \ \ \ \ \ \ \ \ \ \ }+\left.  \left.
{\displaystyle\int_{Z}}
\left(  r\left(  s,z;t\right)  ^{\top}F\left(  z\right)  v+\frac{1}{2}v^{\top
}F\left(  z\right)  ^{\top}P\left(  s;t\right)  F\left(  z\right)  v\right)
\theta\left(  dz\right)  \right\}  1_{\left[  t,t+\varepsilon\right)  }\left(
s\right)  ds\right]  +o\left(  \varepsilon\right)  ,
\end{align*}
which is equivalent to $\left(  3.19\right)  .$\eop

\textbf{A.2.} \textbf{Proof of Theorem 5.2.} For a constant $\beta>0,$ to be
fixed later, we introduce the following norm$,$ for $f\left(  .,.\right)  \in
C\left(  D\left[  0,T\right]  ,%
\mathbb{R}
\right)  $%
\[
\left\Vert f\right\Vert _{\infty,\beta}=\sup_{\left(  t,s\right)  \in D\left[
0,T\right]  }\left\vert e^{-\beta\left(  T-s\right)  }f\left(  t,s\right)
\right\vert ,
\]
it is easy to check that $e^{-\beta T}\left\Vert f\right\Vert _{\infty}%
\leq\left\Vert f\right\Vert _{\infty,\beta}\leq\left\Vert f\right\Vert
_{\infty},$ for every $f\in C\left(  D\left[  0,T\right]  ,%
\mathbb{R}
\right)  $, hence the norm $\left\Vert .\right\Vert _{\infty,\beta}$ is
equivalent to $\left\Vert .\right\Vert _{\infty}$\textbf{ }on the Banach space
$C\left(  D\left[  0,T\right]  ,%
\mathbb{R}
\right)  .$We introduce the following nonlinear operator, $\tilde{L}\left[
.\right]  :C\left(  D\left[  0,T\right]  ,%
\mathbb{R}
^{+}\right)  \rightarrow C\left(  D\left[  0,T\right]  ,%
\mathbb{R}
^{+}\right)  ,$ such that for all $f\left(  .,.\right)  \in C\left(  D\left[
0,T\right]  ,%
\mathbb{R}
^{+}\right)  ,$ we have%
\[
\tilde{L}\left[  f\left(  t,s\right)  \right]  =\tilde{L}\left[  f\right]
\left(  t,s\right)  =h\left(  t\right)  e^{%
{\textstyle\int_{s}^{T}}
\left(  2a-b^{2}\Upsilon\left(  r,r\right)  -b^{2}f\left(  r,r\right)
\right)  dr}.
\]

Then our problem is equivalent to a fixed point problem for the operator
$\tilde{L}\left[  .\right]  $ in the closed subset $C\left(  D\left[
0,T\right]  ,%
\mathbb{R}
^{+}\right)  $ of the Banach space $\left(  C\left(  D\left[  0,T\right]  ,%
\mathbb{R}
\right)  ,\left\Vert .\right\Vert _{\infty,\beta}\right)  .$

1) \textbf{Existence of solution. }It is clear that $\tilde{L}\left[
.\right]  $ is well defined.\textbf{ }Now, consider $f_{1},f_{2}\in C\left(
\left[  0,T\right]  ,%
\mathbb{R}
^{+}\right)  $, then%
\begin{equation}
\tilde{L}\left[  f_{1}\right]  \left(  t,s\right)  -\tilde{L}\left[
f_{2}\right]  \left(  t,s\right)  =h\left(  t\right)  e^{%
{\textstyle\int_{s}^{T}}
\left(  2a-b^{2}\Upsilon\left(  \tau,\tau\right)  \right)  d\tau}\left(
e^{-b^{2}%
{\textstyle\int_{t}^{T}}
f_{1}\left(  \tau,\tau\right)  d\tau}-e^{-b^{2}%
{\textstyle\int_{t}^{T}}
f_{2}\left(  \tau,\tau\right)  d\tau}\right)  , \tag{6.8}%
\end{equation}
we put $\lambda\left(  t,s\right)  =h\left(  t\right)  e^{%
{\textstyle\int_{s}^{T}}
\left(  2a-b^{2}\Upsilon\left(  \tau,\tau\right)  \right)  d\tau},\forall$
$\left(  t,s\right)  \in D\left[  0,T\right]  ,$ obviousely $\lambda\left(
.,.\right)  $ is uniformely bounded. Then there exists some constant $K>0,$
such that
\[
\left\vert \tilde{L}\left[  f_{1}\right]  \left(  t,s\right)  -\tilde
{L}\left[  f_{2}\right]  \left(  t,s\right)  \right\vert \leq K\left\vert
e^{-b^{2}%
{\textstyle\int_{s}^{T}}
f_{1}\left(  \tau,\tau\right)  d\tau}-e^{-b^{2}%
{\textstyle\int_{s}^{T}}
f_{2}\left(  \tau,\tau\right)  d\tau}\right\vert ,
\]
moreover, since $\left\vert e^{-x}-e^{-y}\right\vert \leq\left\vert
x-y\right\vert ,$ $\forall x,y\in\left[  0,\infty\right)  $, then%
\begin{equation}
\left\vert \tilde{L}\left[  f_{1}\right]  \left(  t,s\right)  -\tilde
{L}\left[  f_{2}\right]  \left(  t,s\right)  \right\vert \leq Kb^{2}%
{\displaystyle\int_{s}^{T}}
\left\vert f_{1}\left(  \tau,\tau\right)  -f_{2}\left(  \tau,\tau\right)
\right\vert d\tau, \tag{6.9}%
\end{equation}
thus%
\begin{align*}
e^{-\beta\left(  T-s\right)  }\left\vert \tilde{L}\left[  f_{1}\right]
\left(  t,s\right)  -\tilde{L}\left[  f_{2}\right]  \left(  t,s\right)
\right\vert  &  \leq e^{-\beta\left(  T-s\right)  }Kb^{2}%
{\displaystyle\int_{s}^{T}}
\left\vert f_{1}\left(  \tau,\tau\right)  -f_{2}\left(  \tau,\tau\right)
\right\vert d\tau,\\
&  =e^{-\beta\left(  T-s\right)  }Kb^{2}%
{\displaystyle\int_{s}^{T}}
e^{\beta\left(  T-\tau\right)  }e^{-\beta\left(  T-\tau\right)  }\left\vert
f_{1}\left(  \tau,\tau\right)  -f_{2}\left(  \tau,\tau\right)  \right\vert
d\tau,\\
&  \leq\frac{Kb^{2}\left(  1-e^{-\beta\left(  T-s\right)  }\right)  }{\beta
}\left\Vert f_{1}-f_{2}\right\Vert _{\infty,\beta},
\end{align*}
hence
\[
\left\Vert \tilde{L}\left[  f_{1}\right]  -\tilde{L}\left[  f_{2}\right]
\right\Vert _{\infty,\beta}\leq\frac{Kb^{2}\left(  1-e^{-\beta T}\right)
}{\beta}\left\Vert f_{1}-f_{2}\right\Vert _{\infty,\beta}.
\]

Therefore $\tilde{L}\left[  .\right]  $ is a contraction mapping for $\beta$
large enough.

2) \textbf{Uniqueness of solution. }Let $f_{1},f_{2}\in C\left(  D\left[
0,T\right]  ,%
\mathbb{R}
^{+}\right)  $ be two solutions$,$ then%
\[
f_{1}\left(  t,s\right)  =\tilde{L}\left[  f_{1}\right]  \left(  t,s\right)
\text{ and }f_{2}\left(  t,s\right)  =\tilde{L}\left[  f_{2}\right]  \left(
t,s\right)  ,\text{ }\forall\left(  t,s\right)  \in D\left[  0,T\right]  .
\]

From $\left(  6.9\right)  $ we have%
\[
\left\vert f_{1}\left(  s,s\right)  -f_{2}\left(  s,s\right)  \right\vert \leq
Kb^{2}%
{\displaystyle\int_{s}^{T}}
\left\vert f_{1}\left(  \tau,\tau\right)  -f_{2}\left(  \tau,\tau\right)
\right\vert d\tau,\forall s\in\left[  0,T\right]  ,
\]
therefore, by Gronwall Lemma, we conclude that $\left\vert f_{1}\left(
s,s\right)  -f_{2}\left(  s,s\right)  \right\vert =0,\forall s\in\left[
0,T\right]  ,$ hence%
\begin{align*}
f_{1}\left(  t,s\right)   &  =\tilde{L}\left[  f_{1}\right]  \left(
t,s\right) \\
&  =h\left(  t\right)  e^{%
{\displaystyle\int_{s}^{T}}
\left(  2a-b^{2}\Upsilon\left(  \tau,\tau\right)  -b^{2}f_{1}\left(  \tau
,\tau\right)  \right)  d\tau}\\
&  =h\left(  t\right)  e^{%
{\displaystyle\int_{s}^{T}}
\left(  2a-b^{2}\Upsilon\left(  \tau,\tau\right)  -b^{2}f_{2}\left(  \tau
,\tau\right)  \right)  d\tau}\\
&  =f_{2}\left(  t,s\right)  ,\text{ }\forall\left(  t,s\right)  \in D\left[
0,T\right]
\end{align*}

This completes the proof.\eop


\begin{thebibliography}{99}                                                                                               %


\bibitem {21}S. Basak and G. Chabakauri, \textit{Dynamic mean-variance asset
allocation}, Review of Financial Studies, \textbf{23 }(2010), 2970-3016.

\bibitem {8}T. Bj\"{o}rk and A. Murgoci, \textit{A general theory of Markovian
time-inconsistent stochastic control problems}, SSRN:1694759, (2008).

\bibitem {13}T. Bj\"{o}rk, A. Murgoci and X.Y. Zhou, \textit{Mean-variance
portfolio optimization with state-dependent risk aversion}, Mathematical
Finance, \textbf{24}(1) (2014), 1-24.

\bibitem {19}C. Czichowsky, \textit{Time-consistent mean-variance porftolio
selection in discrete and continuous time}, arXiv:1205.4748v1, (2012).

\bibitem {14}B. Djehiche and M. Huang, \textit{A characterization of sub-game
perfect Nash equilibria for SDEs of mean field type}, Dynamic Games and
Applications, DOI:10.1007/s13235-015-0140-8.

\bibitem {1}I. Ekeland and A. Lazrak, \textit{Equilibrium policies when
preferences are time-inconsistent}, arXiv:0808.3790v1, (2008).

\bibitem {2}I. Ekeland and T.A. Pirvu, \textit{Investment and consumption
without commitment,} Mathematics and Financial Economics, 2 (2008), 57-86.

\bibitem {11}S. M. Goldman, \textit{Consistent plans}, Rev. Financial Stud.,
47 (1980), 533-537.

\bibitem {18}Y. Hu, H. Jin and X.Y. Zhou, \textit{Time-inconsistent stochastic
linear quadratic control}. SIAM J. Control Optim., \textbf{50}(3) (2012), 1548-1572.

\bibitem {5}P. Krusell and A. Smith, \textit{Consumption and savings decisions
with quasi-geometric discounting}, Econometrica, \textbf{71} (2003), 366--375.

\bibitem {23}Q. Meng, \textit{General linear quadratic optimal stochastic
control problem driven by a Brownian motion and a Poisson random martingale
measure with random coefficients}, Stochastic Analysis and Applications, DOI: 10.1080/07362994.2013.845106

\bibitem {12}B. $\varnothing$ksendal and A. Sulem, Applied Stochastic Control
of Jump Diffusions, Second Edition, Springer, $(2007).$

\bibitem {17}S. Peng, \textit{A general stochastic maximum principle for
optimal control problems,} SIAM Journal on Control and Optimization,
\textbf{28} (1990), 966--979.

\bibitem {10}E.S. Phelps and R.A. Pollak, \textit{On second-best national
saving and game-equilibrium growth.} Review of Economic Studies, \textbf{35}
(1968), 185-199.

\bibitem {7}R. Pollak, \textit{Consistent planning}, Rev. Financial Stud.,
\textbf{35} (1968), 185--199.

\bibitem {6}R. Strotz, \textit{Myopia and inconsistency in dynamic utility
maximization}, Review of Economic Studies, \textbf{23} (1955), 165--180.

\bibitem {15}S. Tang and X. Li, \textit{Necessary conditions for optimal
control for stochastic systems with random jumps}. SIAM J. Control Optim.,
\textbf{32} (1994) 1447-1475.

\bibitem {22}Z. Wu and X. Wang, \textit{FBSDE with Poisson process and its
application to linear quadratic stochastic optimal control problem with random
jumps}, Acta Automatica Sinica, \textbf{29} (2003), 821--826.

\bibitem {3}J. Yong, \textit{A deterministic linear quadratic
time-inconsistent optimal control problem,} Math. Control Related Fields,
\textbf{1} (2011), 83-118.

\bibitem {20}J. Yong, \textit{linear quadratic optimal control problems for
mean-field stochastic differential equations: Time-consistent solutions,}
arXiv:1304.3964, (2013).

\bibitem {4}J. Yong, \textit{Time-inconsistent optimal control problems and
the equilibrium HJB equation,} Math. Control Related Fields, \textbf{2}(3)
(2012), 271-329.

\bibitem {9}J. Yong and X.Y. Zhou, Stochastic Controls: Hamiltonian Systems
and HJB Equations. Springer-Verlag, New York. (1999)

\bibitem {16}X.Y. Zhou and D. Li, \textit{Continuous-Time Mean-Variance
Portfolio Selection}: \textit{A Stochastic LQ Framework,} Appl. Math. Optim.,
\textbf{42} (2000), 19-33.
\end{thebibliography}
\end{document}